\documentclass{article}
\usepackage{amsmath,amsthm,amsfonts,amssymb}
\usepackage[dvips]{graphicx}
\usepackage[english,french]{babel}

\newtheorem{e-proposition}[theorem]{Proposition}

\newtheorem{e-definition}[theorem]{Definition\rm}

\def\og{\leavevmode\raise.3ex\hbox{$\scriptscriptstyle\langle\!\langle$~}}
\def\fg{\leavevmode\raise.3ex\hbox{~$\!\scriptscriptstyle\,\rangle\!\rangle$}}

\begin{document}
%\begin{frontmatter}
\selectlanguage{english}
\begin{center}
\textbf{Reproductive strong solutions of Navier-Stokes equations
with non homogeneous boundary conditions}
\end{center}
\begin{center}
\textbf{ Ch\'erif Amrouche $^{1} $},
%\ead{cherif.amrouche@univ-pau.fr}
\textbf{Macaire Batchi $^{1,2,3} $},
%\ead{ma.batchi@etud.univ-pau.fr}
\textbf{Jean Batina $^{2} $}.
%\ead{jean.batina@univ-pau.fr}
\end{center}
\begin{center}
\textit{1\quad Laboratoire de Math\'ematiques Appliqu\'ees CNRS UMR
5142}
\\
\textit{2\quad Laboratoire de Thermique Energ\'etique et
Proc\'ed\'es}
\\
\textit{Universit\'e de Pau et des Pays de l'Adour}\\
\textit{Avenue de l'Universit\'e 64000 Pau, France}\\
\textit{3\quad Universit\'e Marien NGouabi-Facult\'e des Sciences}\\
\textit{B.P.:69 Brazzaville, Congo}

\end{center}
\begin{abstract}
\selectlanguage{english} The object of the present paper is to show
the existence and the uniqueness of a reproductive strong solution
of the Navier-Stokes equations, i.e. the solution $\boldsymbol{u} $
belongs to
$\text{}\mathbf{L}%
^{\infty }\left( 0,T;V\right)\, \cap \,\mathbf{L}^{2}\left( 0,T;\mathbf{H}%
^{2}\left( \Omega \right) \right)$  and
satisfies the property $\boldsymbol{u}%
\left( \boldsymbol{x,\,}T\right) =\boldsymbol{u}%
\left( \boldsymbol{x,\,}0\right) =\boldsymbol{u}_{0}\left(
\boldsymbol{x}\right)$. One considers the case of an incompressible
fluid in two dimensions with nonhomogeneous boundary conditions, and
external forces are neglected.
\end{abstract}
%\end{frontmatter}
%
\smallskip
\noindent\textbf{Key Words}: Navier-Stokes equations, incompressible
fluid, reproductive solution, nonhomogeneous boundary conditions.
\\\\
\medskip
\noindent\textbf{Mathematics Subject %
Classification (2000)}: 35K, 76D03, 76D03 

%\section{Le probl\`{e}me dynamique}

\section{Introduction and notations}

Let $\Omega $ be an open and bounded domain of $%
%TCIMACRO{\U{211d} }%
%BeginExpansion
\mathbb{R}
%EndExpansion
^{2},$ with a sufficiently smooth boundary $\Gamma $; and let us consider the Navier-Stokes equations:%
\begin{equation}
\left\{
\begin{array}{lll}
\dfrac{\partial \boldsymbol{v}}{\partial t}-\nu \triangle \boldsymbol{v}%
\text{ }\mathbf{+}\text{ }\boldsymbol{v}\mathbf{.\nabla
}\boldsymbol{v}\text{ }\mathbf{+}\text{ }\mathbf{\nabla }p=0\quad &
\text{in\quad } &
Q_{T}=\Omega \times \left] 0,T\right[ , \\
\text{div }\boldsymbol{v}\text{ }\mathbf{=}\text{ }0 & \text{in} &
Q_{T},
\\
\boldsymbol{v}\text{ }\mathbf{=}\text{ }\boldsymbol{g} & \text{on} &
\Sigma
_{T}=\Gamma \times \left] 0,T\right[ , \\
\boldsymbol{v}\mathbf{(}0\mathbf{)=}\text{ }\boldsymbol{v}_{0} &
\text{in}
& \Omega .%
\end{array}%
\right.
\end{equation}

\bigskip

\noindent where $\boldsymbol{g}$ , $\boldsymbol{v}_{0}$ and $T>0$ are given.%
We suppose that :%
\begin{equation}
\begin{array}{ll}
\text{div }\boldsymbol{v}_{0}\text{ }\mathbf{=}\text{ }0\text{ }\
\text{in%
\quad }\Omega , & \boldsymbol{v}_{0}.\boldsymbol{n}=0\text{ \quad on \quad }%
\Gamma ,%
\end{array}
\end{equation}%
and
\begin{equation}
\boldsymbol{g}.\boldsymbol{n}=0\text{ \quad on \quad }\Sigma _{T}.
\end{equation}

\bigskip
One is interested on one hand by the existence of strong solutions
of system (1).  On the other hand, one seeks data conditions to
establish the existence of a reproductive solution generalizing the
concept of a periodic solution.  Kaniel and Shinbrot
$%
\left[ 5\right] $ showed the existence of these solutions for system
(1) in dimensions 2 and 3 with external forces but zero boundary
condition i.e. $\boldsymbol{g}=0.$ With another approach using
semigroups, one
can also point out the work of Takeshita $\left[ 10%
\right] $ in dimension 2.

\smallskip

We need to introduce the following functional spaces, with $r$ and
$s$ positive numbers:

\begin{equation*}
\mathbf{H}^{r,s}(Q_{T})\mathit{\ }=\mathbf{L}^{2}\left( \left] 0,T\right[ ;%
\mathbf{H}^{r}(\Omega )\right) \cap \mathbf{H}^{s}\left( \left] 0,T\right[ ;%
\mathbf{L}^{2}(\Omega )\right)
\end{equation*}

\bigskip

\noindent These are Hilbert spaces for the norm

\begin{equation*}
\left\Vert \boldsymbol{v}\right\Vert _{\mathbf{H}^{r,s}(Q_{T})\mathit{\ }%
}=\left( \int\limits_{0}^{T}\left\Vert \boldsymbol{v}\mathbf{(}t\mathbf{)}%
\right\Vert _{\mathbf{H}^{r}(\Omega )}^{2}dt+\left\Vert \boldsymbol{v}%
\right\Vert _{\mathbf{H}^{s}\left( \left] 0,T\right[ ;\mathbf{L}^{2}(\Omega
)\right) }^{2}\right) ^{1/2}.
\end{equation*}

\bigskip

\noindent Let us recall that for $s=1$, for example,

\begin{equation*}
\left\Vert \boldsymbol{v}\right\Vert _{\mathbf{H}^{1}\left( \left] 0,T\right[
;\mathbf{L}^{2}(\Omega )\right) \mathit{\ }}=\left[ \int\limits_{0}^{T}%
\left( \left\Vert \boldsymbol{v}\mathbf{(}t\mathbf{)}\right\Vert _{\mathbf{L}%
^{2}(\Omega )}^{2}+\left\Vert \frac{\partial \boldsymbol{v}}{\partial t}%
\right\Vert _{\mathbf{L}^{2}(\Omega )}^{2}\right) dt\right] ^{1/2}.
\end{equation*}

\bigskip

\noindent In the same manner one defines spaces $\mathbf{H}%
^{r,s}(\Sigma _{T}).\mathit{\ \bigskip }$

\noindent We now introduce the following spaces:%
\begin{equation*}
\begin{array}{l}
\mathcal{V}=\left\{ \boldsymbol{v}\in \mathcal{D}(\Omega )^{2};\text{ div }%
\boldsymbol{v}=0\text{ in \quad }\Omega \right\} , \\
\mathrm{H}=\left\{ \boldsymbol{v}\in \mathbf{L}^{2}(\Omega );\text{ div }%
\boldsymbol{v}=0\text{ in \quad }\Omega ,\ \boldsymbol{v}\mathbf{.}%
\boldsymbol{n}=0\text{ on \quad }\Gamma \right\} , \\
V=\left\{ \boldsymbol{v}\in \mathbf{H}_{0}^{1}(\Omega );\text{ div }%
\boldsymbol{v}=0\text{ in \quad }\Omega \right\} ,%
\end{array}%
\end{equation*}

\smallskip
\noindent Let us recall that $\mathcal{V}$ is dense in $\mathrm{H}$
and $V$ for their respective topologies.

\smallskip

Here, $\mathcal{D}(\Omega )\mathcal{\ }$is the class of
$\mathcal{C}^{\infty }$  functions with compact support in $\Omega
.$ The
notations $\left( \mathbf{.},\mathbf{.}\right) $ et $\left( \left( \mathbf{.}%
,\mathbf{.}\right) \right) $ indicate the scalar products in $%
\mathbf{L}^{2}(\Omega )$ and in $\mathbf{H}_{0}^{1}(\Omega )$
respectively, and $\left\vert .\right\vert $ et $\left\Vert
.\right\Vert $ the associated norms.\medskip

In the order to solve problem (1), we will have to remove boundary
condition $\boldsymbol{g}\mathbf{.\;}$ and consider a new problem
with zero boundary condition. We note that if $\boldsymbol{v}$
$\mathbf{\in H}^{2,1}(Q_{T})\mathit{\ }$ is solution of (1), then
thanks to the Aubin compactness lemma (see J.L.\
Lions $\left[ 8\right] ,$ R.\ Temam $\left[ 11\right] $ ) one will have $\qquad $%
\begin{equation*}
\boldsymbol{v}\text{ }\mathbf{\in }\text{ }\mathcal{C}^{0}\left( \left[ 0,T%
\right] ;\mathbf{H}^{1}(\Omega )\right) \hookrightarrow \mathcal{C}%
^{0}\left( \left[ 0,T\right] ;\mathbf{H}^{1/2}(\Gamma )\right)
\end{equation*}

\noindent So that a necessary condition for $\boldsymbol{v}$ to
exist is that:

\begin{equation}
\begin{array}{ll}
\boldsymbol{g}\left( \boldsymbol{x}\mathbf{,}0\right) =\boldsymbol{v}%
_{0}\left( \boldsymbol{x}\right) , & \text{$\boldsymbol{x}$}\in \Gamma .%
\end{array}
\end{equation}%
\medskip

\noindent Combining (2)-(4), one has:
\begin{equation*}
\begin{array}{ll}
\boldsymbol{g}.\boldsymbol{n}=0\text{ } & \text{on \ }\Gamma \times
\left[
0,T\right[ .%
\end{array}%
\end{equation*}

\noindent The following lemma allows us to state hypotheses on
$\boldsymbol{g}$ (voir Lions-Magenes $\left[ 7\right] $)$.\bigskip $
\smallskip

\noindent \textbf{Lemma 1.1.} \textit{Suppose that (4) takes place
and let}

\begin{equation}
\boldsymbol{g} \mathbf{\in H}^{3/2,3/4}(\Sigma _{T}),\text{ \quad }%
\boldsymbol{v}_{0}\in \mathbf{H}^{1}(\Omega ).
\end{equation}%
\smallskip
\noindent \textit{Then there exists a function
}$\mathbf{R}$\textbf{\ }$\in \mathbf{H}^{2,1}(Q_{T})$\textit{\ such
that}
\begin{equation}
\mathbf{R=}\boldsymbol{g}\text{ \thinspace on \quad }\Sigma
_{T}\text{ et \ }\mathbf{R}\left( 0\right) \text{ }\mathbf{=}\text{
}\boldsymbol{v}_{0}\text{ in \quad }\Omega ,
\end{equation}%
%\smallskip
\noindent \textit{and satisfying the estimates}%
\begin{equation}
\left\Vert \mathbf{R}\right\Vert _{\mathbf{H}^{2,1}(Q_{T})}\leq C\left(
\left\Vert \boldsymbol{g}\right\Vert _{\mathbf{H}^{3/2,3/4}(\Sigma
_{T})}+\left\Vert \boldsymbol{v}_{0}\right\Vert _{\mathbf{H}^{1}(\Omega )%
\mathit{\ }}\right) .\square
\end{equation}%
\bigskip\ \ \ We now consider the problem:\medskip

\noindent \noindent For a given $\boldsymbol{g}$ verifying (5), one seeks %
$\left( \boldsymbol{u}\mathbf{,}\text{ }q\right) $ which satisfies %
\begin{equation}
\left\{
\begin{array}{lll}
\dfrac{\partial \boldsymbol{u}}{\partial t}-\nu \triangle \boldsymbol{u}%
\text{ }\mathbf{+}\text{ }\mathbf{\nabla }q\text{ }\mathbf{=}\text{ }0\quad
& \text{in \quad } & Q_{T}, \\
\text{div }\boldsymbol{u}\text{ }\mathbf{=}\text{ div }\mathbf{R} & \text{%
in} & Q_{T}, \\
\boldsymbol{u}\text{ }\mathbf{=}\text{ }0 & \text{on} & \Sigma _{T}, \\
\boldsymbol{u}\mathbf{(}0\mathbf{)=0} & \text{in} & \Omega .%
\end{array}%
\right.
\end{equation}%
\bigskip\ The following proposition holds (see  Dautray-Lions $\left[ 2%
\right] ,$ O.\ A.\ Ladyzhenskaya $\left[ 6\right] ,$ V.A. Solonnikov
$\left[ 9\right] $) $:\bigskip $

\noindent \textbf{Proposition 1.2. }\textit{We suppose that
(5)holds,}
\begin{equation}
\begin{array}{l}
\text{div }\boldsymbol{v}_{0}\text{ }\mathbf{=}\text{ }0\text{ }\ \text{%
\textit{on} \thinspace }\Omega ,\text{ \thinspace }\boldsymbol{v}_{0}.%
\boldsymbol{n}=0\text{ \thinspace \textit{in} \thinspace }\Gamma
,\textit{and} \text{  }%
\boldsymbol{g}.\boldsymbol{n}=0\text{ \thinspace%
\textit{in } }\Sigma_{T}.%
\end{array}
\end{equation}%

\medskip \textit{Then problem (8) has an unique solution}
$\left( \boldsymbol{u}\mathbf{,}\text{ }q\right) $\textit{\ such that%
}
\begin{equation*}
\begin{array}{l}
\boldsymbol{u}\in \mathbf{H}^{2,1}(Q_{T}),\text{ \hspace{1cm}}q\in
L^{2}\left( 0,T;H^{1}(\Omega )^{2}\right)%
\end{array}%
\end{equation*}%
\textit{with the estimates}%
\begin{equation}
\begin{array}{l}
\left\Vert \boldsymbol{u}\right\Vert _{\mathbf{H}^{2,1}(Q_{T})}+\left\Vert
q\right\Vert _{L^{2}\left( 0,T;H^{1}(\Omega )^{2}\right) \mathit{\ }}\leq
C\left( \left\Vert \boldsymbol{g}\right\Vert _{\mathbf{H}^{3/2,3/4}(\Sigma
_{T})}+\left\Vert \boldsymbol{v}_{0}\right\Vert _{\mathbf{H}^{1}(\Omega )%
\mathit{\ }}\right) .\square%
\end{array}
\end{equation}%
\medskip
%\bigskip
\noindent Thus the function defined by
\begin{equation}
\mathbf{G=R}-\text{ }\boldsymbol{u}\text{ \qquad \qquad\ \ \
\thinspace \thinspace \thinspace in }Q_{T}
\end{equation}%
%\medskip
\noindent satisfies the estimates (7) and
\begin{equation}
\begin{array}{ll}
\text{div }\mathbf{G=}\text{ }0\text{ \qquad \qquad\ \ \ \thinspace
\thinspace \thinspace \quad } & \text{in }Q_{T},%
\end{array}
\end{equation}%
%\medskip
\begin{equation}
\begin{array}{ll}
\mathbf{G=}\text{ }\boldsymbol{g}\text{ \qquad \qquad\ \ \ \thinspace
\thinspace \thinspace \quad } & \ \ \ \ \ \text{on }\Sigma _{T},%
\end{array}
\end{equation}%
%\medskip \medskip
\begin{equation}
\begin{array}{ll}
\mathbf{G}\left( \boldsymbol{x}\mathbf{,}0\right) \text{ }\mathbf{=}\text{ }%
\boldsymbol{v}\left( \boldsymbol{x}\mathbf{,}0\right) \text{ \qquad\ \ \ \ \
} & \text{$\boldsymbol{x}$ }\in \Omega .%
\end{array}
\end{equation}%
This yields the following lemma: \bigskip

\noindent \textbf{Lemma 1.3. }\textit{Let }$\boldsymbol{g}\mathbf{\ }$ and $%
\boldsymbol{v}_{0}$ \textit{satisfy }(4), (5) and (9). \textit{%
Then there exists }$\mathbf{G\in H}^{2,1}(Q_{T})$ \textit{satisfying
(12)-(14)} \textit{and the estimate}

\begin{equation*}
\left\Vert \mathbf{G}\right\Vert _{\mathbf{H}^{2,1}(Q_{T})}\leq C\left(
\left\Vert \boldsymbol{g}\right\Vert _{\mathbf{H}^{3/2,3/4}(\Sigma
_{T})}+\left\Vert \boldsymbol{v}_{0}\right\Vert _{\mathbf{H}^{1}(\Omega )%
\mathit{\ }}\right) .\square
\end{equation*}%
\bigskip Moreover, one has the next lemma \medskip

\noindent \textbf{Lemma 1.4. }\textit{Let }$\varepsilon >0,$\textit{and let \ }$%
\boldsymbol{g}\mathbf{\ }$ and $\boldsymbol{v}_{0}$ \textit{satisfy
the hypotheses of lemma 1.3. Then there exists }$\mathbf{G}_{\varepsilon }$ $%
\mathbf{\in H}^{2,1}(Q_{T})$ \textit{such that}%
\begin{eqnarray*}
&&%
\begin{array}{ll}
\text{div }\mathbf{G}_{\varepsilon }\text{ }\mathbf{=}\text{ }0\text{ } &
\text{in \quad }Q_{T},%
\end{array}
\\
&& \\
&&%
\begin{array}{ll}
\mathbf{G}_{\varepsilon }\text{ }\mathbf{=}\boldsymbol{g}\text{ \quad\ } &
\text{on \ \ \ \quad }\Sigma _{T},%
\end{array}%
\end{eqnarray*}%
%\textit{and satisfying}
\begin{equation*}
\left\Vert \mathbf{G}_{\varepsilon }\left( .,0\right) \right\Vert _{\mathbf{H%
}^{1}(\Omega )}\leq C_{\varepsilon }\left\Vert \mathbf{G}\left( .,0\right)
\right\Vert _{\mathbf{H}^{1}(\Omega )}
\end{equation*}%
\textit{and }
\begin{equation*}
\forall \boldsymbol{v}\text{ }\mathbf{\in }\text{ }V\mathbf{,\quad }%
\left\vert b\left( \boldsymbol{v}\mathbf{,G}_{\varepsilon }\left( t\right) ,%
\boldsymbol{v}\right) \right\vert \leq \beta (\varepsilon ,t)\left\Vert
\nabla \boldsymbol{v}\right\Vert _{\mathbf{L}^{2}\left( \Omega \right) }^{2}
\end{equation*}%
\textit{with}
\begin{equation*}
\underset{t\in \left[ 0,T\right] }{\sup }\beta (\varepsilon ,t)\rightarrow 0%
\text{ \textit{when } }\varepsilon \rightarrow 0.
\end{equation*}%
\medskip

\noindent \textit{Moreover, there exists an increasing function L : }$%
%TCIMACRO{\U{211d} }%
%BeginExpansion
\mathbb{R}
%EndExpansion
^{+}\rightarrow
%TCIMACRO{\U{211d} }%
%BeginExpansion
\mathbb{R}
%EndExpansion
^{+},$\textit{\ not depending on }$\varepsilon ,$\textit{\ such
that\ }
\begin{equation*}
\left\Vert \mathbf{G}_{\varepsilon }\right\Vert _{\mathbf{H}%
^{2,1}(Q_{T})}\leq L\left( \frac{\varepsilon }{\left\Vert \boldsymbol{g}%
\right\Vert _{\mathbf{H}^{3/2,3/4}(\Sigma _{T})}+\left\Vert \boldsymbol{v}%
_{0}\right\Vert _{\mathbf{H}^{1}(\Omega )\mathit{\ }}}\right) \left(
\left\Vert \boldsymbol{g}\right\Vert _{\mathbf{H}^{3/2,3/4}(\Sigma
_{T})}+\left\Vert \boldsymbol{v}_{0}\right\Vert _{\mathbf{H}^{1}(\Omega
)}\right) .
\end{equation*}%
\bigskip \medskip

\noindent \textbf{Proof.\bigskip }

\noindent \textit{i) Step 1 }: %On pose $\rho \left( x\right) =d\left(
%x,\Gamma \right) .\medskip $
\noindent One takes up again the Hopf construction (see Girault \&
Raviart $\left[
4\right] $, Temam $\left[ 11\right] $, Lions $\left[ 8\right] $, Galdi $%
\left[ 3\right] $ )$.\;$
%On d\'{e}finit d'abord la fonction%
\smallskip\\
\noindent \textit{ii) Step 2 }: The open domain $\Omega $ \ being
smooth,
and since div $\mathbf{G}_{\varepsilon }$ $\mathbf{=}$ $0$ in $Q_{T}$ \ and $%
\mathbf{G}.\boldsymbol{n}=0$ \thinspace on $\Gamma \times \left[
0,T\right[
,$ there exists, for all $t\in \left[ 0,T\right[ $, a function $\psi $ depending on %
 $\boldsymbol{x}$ and $t,$ such that

\begin{equation*}
\mathbf{G}=\mathbf{rot}\text{ }\psi \text{ \quad in }\ \quad \Omega \times %
\left[ 0,T\text{ }\right]
\end{equation*}%
\smallskip

\noindent with $\psi =0$ on $\Gamma \times \left[ 0,T\right[ $,
$\psi \in \mathbf{L}^{2}\left( 0,T;\mathbf{H}^{3}(\Omega )\right) $,
$\dfrac{\partial \psi }{\partial t}\in \mathbf{L}^{2}\left(
0,T;\mathbf{H}^{1}(\Omega )\right) $ and satisfying the estimate

\begin{equation}
\left\Vert \psi \right\Vert _{\mathbf{L}^{2}\left( 0,T;\mathbf{H}^{3}(\Omega
)\right) }+\left\Vert \psi _{t}\right\Vert _{\mathbf{L}^{2}\left( 0,T;%
\mathbf{H}^{1}(\Omega )\right) }\leq C\left\Vert \mathbf{G}\right\Vert _{%
\mathbf{H}^{2,1}(Q_{T})}.
\end{equation}%
%\bigskip

\smallskip

\noindent \textit{iii) Step 3 }: Let%
\begin{equation*}
\mathbf{G}^{\varepsilon }=\mathbf{rot}\left( \theta _{\varepsilon }\text{ }%
\psi \right) .
\end{equation*}%
\bigskip One deduces from the properties  of $\theta _{\varepsilon
}$, for $j=1,2$:

\begin{equation*}
\left\vert \mathbf{G}_{j}^{\varepsilon }(x,t)\right\vert \leq C\left( \dfrac{%
\varepsilon }{\rho \left( x\right) }\left\vert \psi (x,t)\right\vert
+\left\vert \nabla \psi (x,t)\right\vert \right) \qquad \text{if \qquad }%
\rho (x)\leq 2\delta (\varepsilon )
\end{equation*}%
\smallskip

\noindent and $\mathbf{G}_{j}^{\varepsilon }=0$ \quad if $\rho
(x)>2\delta (\varepsilon )$.
\\
\smallskip
\noindent We note that

\begin{equation*}
\psi \in C\left( \left[ 0,T\right] ;\mathbf{H}^{2}(\Omega )\right)
\hookrightarrow C\left( \left[ 0,T\right] ;\mathbf{L}^{\infty }(\Omega
)\right) .
\end{equation*}%
%\smallskip

\noindent Therefore,

\begin{equation*}
\left\vert \mathbf{G}_{j}^{\varepsilon }(x,t)\right\vert \leq C\left( \dfrac{%
\varepsilon }{\rho \left( x\right) }+\left\vert \nabla \psi
(x,t)\right\vert \right) \hspace{1cm}\text{if \qquad }\rho (x)\leq
2\delta (\varepsilon ).
\end{equation*}%
\medskip

\noindent Thus, for all $\boldsymbol{v}\in \mathbf{H}_{0}^{1}(\Omega ),$%
\begin{equation*}
\left\Vert \boldsymbol{v}_{i}\mathbf{G}_{j}^{\varepsilon }\right\Vert _{%
\mathbf{L}^{2}\left( \Omega \right) }\leq C\left[ \varepsilon \left\Vert
\frac{\boldsymbol{v}_{i}}{\rho }\right\Vert _{\mathbf{L}^{2}\left( \Omega
\right) }+\left( \int\limits_{\rho (x)\leq 2\delta (\varepsilon )}%
\boldsymbol{v}_{i}^{2}.\left\vert \nabla \psi \right\vert ^{2}dx\right)
^{1/2}\right]
\end{equation*}

\smallskip

\begin{equation*}
\left\Vert \boldsymbol{v}_{i}\mathbf{G}_{j}^{\varepsilon }\right\Vert _{%
\mathbf{L}^{2}\left( \Omega \right) }\leq C\varepsilon \left\Vert \nabla
\boldsymbol{v}_{i}\right\Vert _{\mathbf{L}^{2}\left( \Omega \right)
}+C\left\Vert \nabla \boldsymbol{v}_{i}\right\Vert _{\mathbf{L}^{2}\left(
\Omega \right) }\times \left( \int\limits_{\rho (x)\leq 2\delta (\varepsilon
)}\left\vert \nabla \psi \right\vert ^{3}dx\right) ^{1/3}
\end{equation*}

\bigskip

\noindent Setting
\begin{equation*}
\beta (\varepsilon ,t)=\left( \int\limits_{\rho (x)\leq 2\delta
(\varepsilon )}\left\vert \nabla \psi \right\vert ^{3}dx\right)
^{1/3},
\end{equation*}

\noindent it's clear that
\begin{equation*}
\underset{\varepsilon \rightarrow 0}{\lim }\beta (\varepsilon ,t)=0
\text{ } uniformly \text{ } on \text{ }\left[ 0,T\right] .
\end{equation*}

\smallskip

\noindent The second inequality of lemma 1.4 is a consequence of %
H\"{o}lder inequality.\ The first inequality follows from %
Hardy inequality for $\mathbf{H}%
_{0}^{1}(\Omega )$-functions and properties of $\theta _{\varepsilon }$ .$%
\square \bigskip $
\section{Existence of strong solutions}

{\Large \smallskip }Let us make a change of the unknown function in
problem (1), by setting

\begin{equation*}
\begin{array}{ll}
\boldsymbol{u}=\boldsymbol{v}-\mathbf{G}_{\varepsilon },\hspace{2.5cm} &
\boldsymbol{u}_{0}=\boldsymbol{v}_{0}-\mathbf{G}_{\varepsilon }\left(
.,0\right) ,%
\end{array}%
\end{equation*}

\bigskip

\noindent where $\mathbf{G}_{\varepsilon }$ is the function given by
lemma 1.4. Problem (1) then becomes:%
\begin{equation}
\left\{
\begin{array}{lll}
\dfrac{\partial \boldsymbol{u}}{\partial t}-\nu \triangle \boldsymbol{u}%
\text{ }\mathbf{+}\text{ }\boldsymbol{u}\mathbf{.\nabla }\boldsymbol{u}\text{
}\mathbf{+}\text{ }\boldsymbol{u}\mathbf{.\nabla \mathbf{G}}_{\varepsilon }%
\text{ }\mathbf{+}\text{ }\mathbf{\mathbf{G}}_{\varepsilon }\mathbf{.\nabla }%
\boldsymbol{u}\text{ }\mathbf{+}\text{ }\mathbf{\mathbf{\nabla }}p\text{ }%
\mathbf{=}\ \boldsymbol{f}_{\varepsilon }\text{ \ } & \text{in\ } &
Q_{T}
\\
\text{div }\boldsymbol{u}\text{ }\mathbf{=}\text{ }0 & \text{in} &
Q_{T}
\\
\boldsymbol{u}\text{ }\mathbf{=}\text{ }0 & \text{on} & \Sigma _{T} \\
\boldsymbol{u}\mathbf{(}0\mathbf{)=}\text{ }\boldsymbol{u}_{0}^{\varepsilon }
& \text{in} & \Omega%
\end{array}%
\right.
\end{equation}

\bigskip

\noindent with%
\begin{equation}
\begin{array}{lll}
\boldsymbol{f}_{\varepsilon }\text{ }\mathbf{=-}\text{ }\dfrac{\partial
\mathbf{G}_{\varepsilon }}{\partial t}+\nu \triangle \mathbf{G}_{\varepsilon
}\text{ }\mathbf{-}\text{ }\mathbf{\mathbf{G}}_{\varepsilon }\mathbf{.\nabla
\mathbf{G}}_{\varepsilon } & \ \text{and } & \boldsymbol{u}_{0}^{\varepsilon }%
\text{ }\mathbf{=}\text{ }\boldsymbol{v}_{0}-\mathbf{G}_{\varepsilon }\left(
.,0\right) .%
\end{array}
\end{equation}

\bigskip

\noindent We note that $\boldsymbol{u}_{0}^{\varepsilon }\in V$ and
\begin{equation}
\left\Vert \boldsymbol{u}_{0}^{\varepsilon }\right\Vert _{\mathbf{H}%
^{1}(\Omega )}\leq C_{\varepsilon }\left( \left\Vert \boldsymbol{g}%
\right\Vert _{\mathbf{H}^{3/2,3/4}(\Sigma _{T})}+\left\Vert \boldsymbol{v}%
_{0}\right\Vert _{\mathbf{H}^{1}(\Omega )\mathit{\ }}\right) .
\end{equation}

\bigskip

\noindent Moreover, $\boldsymbol{f}_{\varepsilon }\in \mathbf{L}%
^{2}\left( 0,T;\mathbf{L}^{2}(\Omega )\right) $ and

\begin{equation}
\left\Vert \boldsymbol{f}_{\varepsilon }\right\Vert _{\mathbf{L}^{2}\left(
0,T;\mathbf{L}^{2}(\Omega )\right) }\leq C_{\varepsilon }\left( \left\Vert
\boldsymbol{g}\right\Vert _{\mathbf{H}^{3/2,3/4}(\Sigma _{T})}+\left\Vert
\boldsymbol{v}_{0}\right\Vert _{\mathbf{H}^{1}(\Omega )\mathit{\ }}\right) .
\end{equation}

\noindent Now we are able to announce and to establish the following
theorem :\bigskip

\noindent \textbf{Theorem 2.1. }\textit{Let \ }$\boldsymbol{v}_{0}$
and $\boldsymbol{g}\mathbf{\ }$ \textit{satisfy the hypotheses of
lemma 1.3. Then problem (16) has a unique solution }$%
\left( \boldsymbol{u}\mathbf{,}\text{ }p\right) $ \textit{such that }%
\begin{equation*}
\begin{array}{lll}
\boldsymbol{u}\mathbf{\in }\text{ }\mathbf{L}^{2}\left( 0,T;\mathbf{H}%
^{2}(\Omega )\right) \cap \mathbf{L}^{\infty }\left( 0,T;V\right) , & \dfrac{%
\partial \boldsymbol{u}}{\partial t}\in \mathbf{L}^{2}\left( 0,T;\mathbf{H}%
\right) , & p\in \mathbf{L}^{2}\left( 0,T;H^{1}(\Omega )\right) ,%
\end{array}%
\end{equation*}

%\bigskip

\noindent $p$ \textit{being unique up to an } $%
\mathbf{L}^{2}\left( 0,T\right)$\textit{-function of the single %
variable t}.\smallskip\medskip

\noindent \textbf{Proof. }

%\smallskip

\subsection{\textit{Approximate solutions }}

\smallskip

\noindent We use the Galerkin method. Let $m\in
%TCIMACRO{\U{2115} }%
%BeginExpansion
\mathbb{N}
%EndExpansion
^{\ast }$ and $\boldsymbol{u}_{0m } \in \left\langle \boldsymbol{w}_{1,}%
\boldsymbol{w}_{2}\mathbf{,...,}\boldsymbol{w}_{m}\right\rangle $
such that
\begin{equation*}
\boldsymbol{u}_{0m }\rightarrow \boldsymbol{u}_{0}^{\varepsilon
}\text{ in  \ }V,\text{ if \ }m\rightarrow \infty ,
\end{equation*}

\bigskip

\noindent where $\boldsymbol{w}_{j}$ are the Stokes operator eigenfunctions . %
For each $m$, one defines an approximate solution of %
(16) by :

\begin{equation}
\left\{
\begin{array}{c}
\boldsymbol{u}_{m}(t)=\sum\limits_{j=1}^{m}g_{jm}(t)\boldsymbol{w}_{j} \\
\left( \boldsymbol{u}_{m}^{\prime }\left( t\right) ,\boldsymbol{w}%
_{j}\right) +\nu \left( \left( \boldsymbol{u}_{m}\left( t\right) \mathbf{,}%
\boldsymbol{w}_{j}\right) \right) +b\left( \boldsymbol{u}_{m}\left( t\right)
\mathbf{,}\boldsymbol{u}_{m}\left( t\right) \mathbf{,}\boldsymbol{w}%
_{j}\right) \\
+b\left( \boldsymbol{u}_{m}\left( t\right) \mathbf{,\mathbf{G}}_{\varepsilon
}\left( t\right) \mathbf{,}\boldsymbol{w}_{j}\right) +b\left( \mathbf{G}%
_{\varepsilon }\left( t\right) \mathbf{,}\boldsymbol{u}_{m}\left( t\right)
\mathbf{,}\boldsymbol{w}_{j}\right) =\left( \boldsymbol{f}_{\varepsilon }%
\mathbf{\ }\left( t\right) \mathbf{,}\boldsymbol{w}_{j}\right) \\
\boldsymbol{u}_{m}(0)=\boldsymbol{u}_{0m},\text{ \thinspace }j=1,...,m%
\end{array}%
\right.
\end{equation}

\bigskip

\noindent This is a nonlinear differential system  %
of m equations in m unknowns $g_{jm},$ $j=1,...,m:\medskip $

$\sum_{i=1}^{m}\left( \boldsymbol{w}_{i},\boldsymbol{w}_{j}\right)
g_{im}^{\prime }\left( t\right) +\nu \sum_{i=1}^{m}\left( \left(
\boldsymbol{w}_{i},\boldsymbol{w}_{j}\right) \right) g_{im}\left(
t\right)
+\sum_{i,l=1}^{m}b\left( \boldsymbol{w}_{i},\boldsymbol{w}_{l}%
\mathbf{,}\boldsymbol{w}_{j}\right) g_{im}\left( t\right) g_{lm}\left(
t\right) +$

$+\sum_{i=1}^{m}\left[ b\left( \boldsymbol{w}_{i}\mathbf{,\mathbf{G}}%
_{\varepsilon }\left( t\right) \mathbf{,}\boldsymbol{w}_{j}\right)
g_{im}\left( t\right) +b\left( \mathbf{G}_{\varepsilon }\left( t\right)
\mathbf{,}\boldsymbol{w}_{i}\mathbf{,}\boldsymbol{w}_{j}\right) g_{im}\left(
t\right) \right] =\left( \boldsymbol{f}_{\varepsilon }\mathbf{\ }\left(
t\right) \mathbf{,}\boldsymbol{w}_{j}\right) ,$ $j=1,...,m\medskip $

\subsection{\textit{ Estimates I}}

\smallskip

\noindent Let us multiply (20) by $g_{jm}(t)$ and sum over %
$j:\medskip $

$\qquad
\begin{array}{ll}
\dfrac{1}{2}\dfrac{d}{dt}\left\vert \boldsymbol{u}_{m}\left( t\right)
\right\vert ^{2}+\nu \left\Vert \boldsymbol{u}_{m}\left( t\right)
\right\Vert ^{2} & =-b\left( \boldsymbol{u}_{m}\left( t\right) \mathbf{,G}%
_{\varepsilon }\left( t\right) \mathbf{,}\boldsymbol{u}_{m}\left( t\right)
\right) +\left( \boldsymbol{f}_{\varepsilon }\mathbf{\ }\left( t\right)
\mathbf{,}\boldsymbol{u}_{m}\left( t\right) \right) \\
& \leq \left\vert \boldsymbol{f}_{\varepsilon }\mathbf{\ }\left( t\right)
\right\vert \left\Vert \boldsymbol{u}_{m}\left( t\right) \right\Vert
+\left\vert b\left( \boldsymbol{u}_{m}\left( t\right) \mathbf{,G}%
_{\varepsilon }\left( t\right) \mathbf{,}\boldsymbol{u}_{m}\left( t\right)
\right) \right\vert%
\end{array}%
$

$\hspace{1.5cm}\hspace{1.5cm}$ $\ \ \ \ \ \ \ \ \ \ \ \ \ \ \ \ \ \ \ \ \ \
\ \ \ \ \ \ \ \ \ \ \ \ \ \ \ \ \ \ $

\noindent One deduces from lemma 1.4 that :\medskip

$\hspace{1.5cm}\dfrac{1}{2}\dfrac{d}{dt}\left\vert \boldsymbol{u}_{m}\left(
t\right) \right\vert ^{2}+\dfrac{\nu }{2}\left\Vert \boldsymbol{u}_{m}\left(
t\right) \right\Vert ^{2}\leq \dfrac{1}{2\nu C^{2}\left( \Omega \right) }%
\left\vert \boldsymbol{f}_{\varepsilon }\mathbf{\ }\left( t\right)
\right\vert ^{2}+\beta (\varepsilon ,t)\left\Vert \boldsymbol{u}_{m}\left(
t\right) \right\Vert ^{2}.\medskip $

\noindent As $\underset{t\in \left[ 0,T\right] }{\sup }\beta
(\varepsilon ,t)\rightarrow 0$ when $\varepsilon \rightarrow 0,$ for
a fixed and small $\varepsilon >0$, one has:
\begin{equation}
\dfrac{d}{dt}\left\vert \boldsymbol{u}_{m}\left( t\right) \right\vert ^{2}+%
\dfrac{\nu }{2}\left\Vert \boldsymbol{u}_{m}\left( t\right) \right\Vert
^{2}\leq \dfrac{1}{\nu C^{2}\left( \Omega \right) }\left\vert \boldsymbol{f}%
_{\varepsilon }\mathbf{\ }\left( t\right) \right\vert ^{2}.
\end{equation}

\noindent

\bigskip

\noindent Integrating (21) from 0 to s, one deduces that:\medskip

$\qquad
\begin{array}{ll}
\left\vert \boldsymbol{u}_{m}(s)\right\vert ^{2} & \leq \left\vert
\boldsymbol{u}_{0m}\right\vert ^{2}+\dfrac{1}{\nu C^{2}\left( \Omega \right)
}\int_{0}^{s}\left\vert \boldsymbol{f}_{\varepsilon }\mathbf{\ }%
\left( t\right) \right\vert ^{2}dt \\
& \leq \left\vert \boldsymbol{u}_{0}^{\varepsilon }\right\vert ^{2}+\dfrac{1%
}{\nu C^{2}\left( \Omega \right) }\left\Vert \boldsymbol{f}_{\varepsilon }%
\mathbf{\ }\left( t\right) \right\Vert _{\mathbf{L}^{2}\left( 0,T;\mathbf{L}%
^{2}(\Omega )\right) }^{2} \\
& \leq C_{\varepsilon }\left( \left\Vert \boldsymbol{g}\right\Vert _{\mathbf{%
H}^{3/2,3/4}(\Sigma _{T})}^{2}+\left\Vert \boldsymbol{v}_{0}\right\Vert _{%
\mathbf{H}^{1}(\Omega )\mathit{\ }}^{2}\right)%
\end{array}%
\medskip $

\noindent according to (18) and (20).\ Therefore
\begin{equation}
\boldsymbol{u}_{m}\text{ \quad }\in \text{\ }\mathbf{L}%
^{\infty }(0,T;\mathbf{H}),
\end{equation}
\noindent and $\left\{ \boldsymbol{u}_{m}\right\}$ is an equibounded
sequence in $\mathbf{L}^{\infty }(0,T;\mathbf{H})$.
\bigskip

\noindent Next, thanks to (21), one has:
\begin{equation}
\boldsymbol{u}_{m}\text{ \quad }\in \text{\ }\mathbf{L}%
^{2}(0,T;V),
\end{equation}%
\noindent and the sequence $\left\{ \boldsymbol{u}_{m}\right\}$ is
equibounded in $\mathbf{L}^{2}(0,T;\mathbf{V})$.

%\bigskip

\smallskip

\subsection{\textit{Estimates II}}

\smallskip

\noindent Let us multiply (20) by $\lambda _{j}g_{jm}(t)$ and sum
over $j$ :
\begin{equation}
\begin{array}{l}
\dfrac{1}{2}\dfrac{d}{dt}\left\Vert \boldsymbol{u}_{m}\left( t\right)
\right\Vert ^{2}+\nu \left\vert A\boldsymbol{u}_{m}\left( t\right)
\right\vert ^{2}+b\left( \boldsymbol{u}_{m}\left( t\right) \mathbf{,}%
\boldsymbol{u}_{m}\left( t\right) \mathbf{,}\text{ }A\boldsymbol{u}%
_{m}\left( t\right) \right) + \\
b\left( \mathbf{G}_{\varepsilon }\left( t\right) \mathbf{,}\boldsymbol{u}%
_{m}\left( t\right) \mathbf{,}\text{ }A\boldsymbol{u}_{m}\left( t\right)
\right) +b\left( \boldsymbol{u}_{m}\left( t\right) \mathbf{,G}_{\varepsilon
}\left( t\right) \mathbf{,}\text{ }A\boldsymbol{u}_{m}\left( t\right)
\right) =\left( \boldsymbol{f}_{\varepsilon }\mathbf{,}\text{ }A\boldsymbol{u%
}_{m}\left( t\right) \right)%
\end{array}
\end{equation}

\bigskip

\noindent where $A$ is the Stokes operator. Let us begin by
considering the nonlinear terms.

\smallskip

\noindent For the first term, thanks to the Gagliardo-Nirenberg
inequality one has \medskip

$\qquad
\begin{array}{ll}
\left\vert b\left( \boldsymbol{u}_{m}\left( t\right) \mathbf{,}\boldsymbol{u}%
_{m}\left( t\right) ,\text{\textbf{\ }}A\boldsymbol{u}_{m}\left( t\right)
\right) \right\vert & \leq \left\Vert \boldsymbol{u}_{m}\left( t\right)
\right\Vert _{\mathbf{L}^{4}\left( \Omega \right) }\left\Vert \nabla
\boldsymbol{u}_{m}\left( t\right) \right\Vert _{\mathbf{L}^{4}\left( \Omega
\right) }\left\vert A\boldsymbol{u}_{m}\left( t\right) \right\vert \\
& \leq C\left\vert \boldsymbol{u}_{m}\left( t\right) \right\vert
^{1/2}\left\Vert \boldsymbol{u}_{m}\left( t\right) \right\Vert \left\vert A%
\boldsymbol{u}_{m}\left( t\right) \right\vert ^{3/2} \\
& \leq C\left\Vert \boldsymbol{u}_{m}\left( t\right) \right\Vert ^{4}+\dfrac{%
\nu }{8}\left\vert A\boldsymbol{u}_{m}\left( t\right) \right\vert ^{2}.%
\end{array}%
$

\smallskip

\noindent In the same way,\medskip

$\qquad
\begin{array}{ll}
\left\vert b\left( \mathbf{G}_{\varepsilon }\left( t\right) \mathbf{,}%
\boldsymbol{u}_{m}\left( t\right) \mathbf{,}\text{ }A\boldsymbol{u}%
_{m}\left( t\right) \right) \right\vert & \leq \left\Vert \mathbf{G}%
_{\varepsilon }\left( t\right) \right\Vert _{\mathbf{L}^{4}\left( \Omega
\right) }\left\Vert \nabla \boldsymbol{u}_{m}\left( t\right) \right\Vert _{%
\mathbf{L}^{4}\left( \Omega \right) }\left\vert A\boldsymbol{u}_{m}\left(
t\right) \right\vert \\
& \leq C\left\Vert \mathbf{G}_{\varepsilon }\left( t\right) \right\Vert _{%
\mathbf{H}^{1}\left( \Omega \right) }\left\Vert \boldsymbol{u}_{m}\left(
t\right) \right\Vert ^{1/2}\left\vert A\boldsymbol{u}_{m}\left( t\right)
\right\vert ^{3/2} \\
& \leq C\left\Vert \mathbf{G}_{\varepsilon }\left( t\right) \right\Vert _{%
\mathbf{H}^{1}\left( \Omega \right) }^{4}\left\Vert \boldsymbol{u}_{m}\left(
t\right) \right\Vert ^{2}+\dfrac{\nu }{8}\left\vert A\boldsymbol{u}%
_{m}\left( t\right) \right\vert ^{2}.%
\end{array}%
$

$\qquad \qquad $\ $\ \ \ \ \ \ \ \ \ \ \ \ \ \ \ \ \ \ \ \ \ \ \ \ \ \ \ \ \
\ \ \ \ \ \ \ \ \ \ \ \ \ \ $

\noindent We remark that, according to lemma 1.4, one %
has:\medskip

$\qquad
\begin{array}{ll}
\left\Vert \mathbf{G}_{\varepsilon }\right\Vert _{\mathbf{L}^{\infty }\left(
0,T;\mathbf{H}^{1}(\Omega )\right) } & \leq C\left( \left\Vert \boldsymbol{g}%
\right\Vert _{\mathbf{H}^{3/2,3/4}(\Sigma _{T})}+\left\Vert \boldsymbol{v}%
_{0}\right\Vert _{\mathbf{H}^{1}(\Omega )}\right) .%
\end{array}%
$

$\qquad $

\noindent So that\medskip

$\qquad
\begin{array}{ll}
\left\vert b\left( \mathbf{G}_{\varepsilon }\left( t\right) \mathbf{,}%
\boldsymbol{u}_{m}\left( t\right) \mathbf{,}\text{ }A\boldsymbol{u}%
_{m}\left( t\right) \right) \right\vert & \leq C\left\Vert \boldsymbol{u}%
_{m}\left( t\right) \right\Vert ^{2}+\dfrac{\nu }{8}\left\vert A\boldsymbol{u%
}_{m}\left( t\right) \right\vert ^{2}.%
\end{array}%
$

$\qquad $

\noindent Finally,\medskip

$\qquad
\begin{array}{ll}
\left\vert b\left( \boldsymbol{u}_{m}\left( t\right) \mathbf{,G}%
_{\varepsilon }\left( t\right) \mathbf{,}\text{ }A\boldsymbol{u}_{m}\left(
t\right) \right) \right\vert & \leq \left\Vert \boldsymbol{u}_{m}\left(
t\right) \right\Vert _{\mathbf{L}^{4}\left( \Omega \right) }\left\Vert
\nabla \mathbf{G}_{\varepsilon }\left( t\right) \right\Vert _{\mathbf{L}%
^{4}\left( \Omega \right) }\left\vert A\boldsymbol{u}_{m}\left( t\right)
\right\vert \\
& \leq C\left\Vert \boldsymbol{u}_{m}\left( t\right) \right\Vert
^{2}\left\Vert \mathbf{G}_{\varepsilon }\left( t\right) \right\Vert _{%
\mathbf{H}^{2}\left( \Omega \right) }^{2}+\dfrac{\nu }{8}\left\vert A%
\boldsymbol{u}_{m}\left( t\right) \right\vert ^{2}.%
\end{array}%
$

$\qquad $

\noindent Hence,\medskip

$\dfrac{d}{dt}\left\Vert \boldsymbol{u}_{m}\left( t\right) \right\Vert
^{2}+\nu \left\vert A\boldsymbol{u}_{m}\left( t\right) \right\vert ^{2}\leq
\dfrac{C}{\nu }\left\vert \boldsymbol{f}_{\varepsilon }\mathbf{\ }\left(
t\right) \right\vert ^{2}+C\left[ \left\Vert \boldsymbol{u}_{m}\left(
t\right) \right\Vert ^{4}+\left\Vert \boldsymbol{u}_{m}\left( t\right)
\right\Vert ^{2}\left( 1+\left\Vert \mathbf{G}_{\varepsilon }\left( t\right)
\right\Vert _{\mathbf{H}^{2}\left( \Omega \right) }^{2}\right) \right] .$

\smallskip

\noindent Let\medskip

$\qquad \sigma _{m}\left( t\right) =C\left[ \left\Vert \boldsymbol{u}%
_{m}\left( t\right) \right\Vert ^{2}+\left( 1+\left\Vert \mathbf{G}%
_{\varepsilon }\left( t\right) \right\Vert _{\mathbf{H}^{2}\left(
\Omega \right) }^{2}\right) \right]. \medskip $

\noindent One knows that\medskip

$\qquad \sigma _{m}\left( t\right) \in %
\mathbf{L}^{1}\left( 0,T\right) ;\medskip $

\noindent so that, according to the Gronwall lemma and (24), one
has:

\begin{equation}
\boldsymbol{u}_{m}\text{ \quad }\in \text{\ }\mathbf{L}%
^{\infty }\left( 0,T;V\right) \cap \mathbf{L}^{2}\left( 0,T;\mathbf{H}%
^{2}\left( \Omega \right) \right),
\end{equation}%
\noindent and $\left\{ \boldsymbol{u}_{m}\right\}$ is an equibounded
sequence in $\mathbf{L}%
^{\infty }\left( 0,T;V\right) \cap \mathbf{L}^{2}\left( 0,T;\mathbf{H}%
^{2}\left( \Omega \right) \right)$.

%\bigskip

\smallskip

\subsection{\textit{Estimates III}}

\smallskip

\noindent Let us multiply (20) by $g_{jm}^{\prime }(t)$ and sum over j %
from 1 to m.\ Then\medskip

$\qquad
\begin{array}{ll}
\left\vert \boldsymbol{u}_{m}^{\prime }\left( t\right) \right\vert ^{2}= &
\nu \left( A\boldsymbol{u}_{m}\left( t\right) ,\boldsymbol{u}_{m}^{\prime
}\left( t\right) \right) -b\left( \boldsymbol{u}_{m}\left( t\right) \mathbf{,%
}\boldsymbol{u}_{m}\left( t\right) \mathbf{,}\boldsymbol{u}_{m}^{\prime
}\left( t\right) \right) \\
& -b\left( \mathbf{G}_{\varepsilon }\left( t\right) \mathbf{,}\boldsymbol{u}%
_{m}\left( t\right) \mathbf{,}\boldsymbol{u}_{m}^{\prime }\left( t\right)
\right) -b\left( \boldsymbol{u}_{m}\left( t\right) \mathbf{,G}_{\varepsilon
}\left( t\right) \mathbf{,}\boldsymbol{u}_{m}^{\prime }\left( t\right)
\right) +\left( \boldsymbol{f}_{\varepsilon }\mathbf{,}\boldsymbol{u}%
_{m}^{\prime }\left( t\right) \right) .%
\end{array}%
$

\smallskip $\qquad \hspace{1.5cm}$

\noindent From this, one deduces that\medskip

$\qquad
\begin{array}{ll}
\left\vert \boldsymbol{u}_{m}^{\prime }\left( t\right) \right\vert ^{2}\leq
& \nu \left\vert A\boldsymbol{u}_{m}\left( t\right) \right\vert \left\vert
\boldsymbol{u}_{m}^{\prime }\left( t\right) \right\vert +C\left\Vert
\boldsymbol{u}_{m}\left( t\right) \right\Vert _{\mathbf{L}^{4}\left( \Omega
\right) }\left\Vert \nabla \boldsymbol{u}_{m}\left( t\right) \right\Vert _{%
\mathbf{L}^{4}\left( \Omega \right) }\left\vert \boldsymbol{u}_{m}^{\prime
}\left( t\right) \right\vert \\
& +\text{ }C\left\Vert \mathbf{G}_{\varepsilon }\left( t\right) \right\Vert
_{\mathbf{L}^{4}\left( \Omega \right) }\left\Vert \nabla \boldsymbol{u}%
_{m}\left( t\right) \right\Vert _{\mathbf{L}^{4}\left( \Omega \right)
}\left\vert \boldsymbol{u}_{m}^{\prime }\left( t\right) \right\vert \\
& +\text{ }C\left\Vert \boldsymbol{u}_{m}\left( t\right) \right\Vert _{%
\mathbf{L}^{4}\left( \Omega \right) }\left\Vert \nabla \mathbf{G}%
_{\varepsilon }\left( t\right) \right\Vert _{\mathbf{L}^{4}\left( \Omega
\right) }\left\vert \boldsymbol{u}_{m}^{\prime }\left( t\right) \right\vert
+\left\vert \boldsymbol{f}_{\varepsilon }\mathbf{\ }\left( t\right)
\right\vert \left\vert \boldsymbol{u}_{m}^{\prime }\left( t\right)
\right\vert%
\end{array}%
$

$\hspace{1.5cm}$

\noindent Using the Gagliardo-Nirenberg inequality,
estimates (25) and (19), and lemma 1.4 giving the estimate of $%
\mathbf{G}_{\varepsilon },$ one deduces that
\begin{equation}
\boldsymbol{u}_{m}^{\prime }\text{ \quad }\in \text{\ }%
\mathbf{L}^{2}\left( 0,T;\mathbf{H}\right),
\end{equation}
\noindent and $\left\{ \boldsymbol{u}_{m}^{\prime }\right\}$ is an
equibounded sequence in $\mathbf{L}%
^{2}\left( 0,T;H\right) $.
\\
\smallskip

\subsection{\textit{Taking the limit.}}

\smallskip

\noindent It is a consequence of the above estimates that the
sequence $\boldsymbol{u}_{m}$ has a subsequence
$\boldsymbol{u}_{m}$, the same notation being used to avoid
unnecessary notation overload:

\begin{align}
\boldsymbol{u}_{m}& \rightharpoonup \boldsymbol{u}\text{ weakly* \qquad in \quad }%
\mathbf{L}^{\infty }\left( 0,T;V\right) \text{,}   \\
\boldsymbol{u}_{m}& \rightharpoonup \boldsymbol{u}\text{ weakly\quad\qquad in\quad }%
\mathbf{L}^{2}\left( 0,T;\mathbf{H}^{2}\left( \Omega \right) \right) \text{,}   \\
\boldsymbol{u}_{m}^{\prime }& \rightharpoonup \boldsymbol{u}^{\prime
}\text{ weakly \qquad in \quad }\mathbf{L}^{2}\left(
0,T;\mathbf{H}\right) \text{ .}
\end{align}

\bigskip

\noindent But we have a compact embedding\medskip

$\hspace{1.5cm}\left\{ \boldsymbol{v}\text{ }\mathbf{\in }\text{ }\mathbf{L}%
^{2}\left( 0,T;\mathbf{H}^{2}\left( \Omega \right) \cap V\right) ,\text{ }%
V^{\prime }\in \mathbf{L}^{2}\left( 0,T;\mathbf{H}\right) \right\} \underset{%
compact}{\hookrightarrow }\mathbf{L}^{2}\left( 0,T;V\right) $

\smallskip

\noindent So that

\begin{equation}
\boldsymbol{u}_{m}\rightarrow \boldsymbol{u}\text{ strongly \hspace{0.5cm}in \quad }%
\mathbf{L}^{2}\left( 0,T;V\right) \text{ and a.e. in }Q_{T}
\end{equation}

\smallskip

\noindent Let $m_{0}$ be fixed and $\boldsymbol{v}\in \left\langle
\boldsymbol{w}_{1,}\boldsymbol{w}_{2}\mathbf{,...,}\boldsymbol{w}%
_{m_{0}}\right\rangle .$ Let m tend towards +$\infty $ in (20).\
Then \medskip

$\qquad
\begin{array}{ll}
\left( \boldsymbol{u}^{\prime }\left( t\right) ,\boldsymbol{v}\right) +\nu
\left( \left( \boldsymbol{u}\left( t\right) \mathbf{,}\boldsymbol{v}\right)
\right) & +b\left( \boldsymbol{u}\left( t\right) \mathbf{,}\boldsymbol{u}%
\left( t\right) \mathbf{,}\boldsymbol{v}\right) +b\left( \boldsymbol{u}%
\left( t\right) \mathbf{,\mathbf{G}}_{\varepsilon }\left( t\right) \mathbf{,}%
\boldsymbol{v}\right) \\
& +b\left( \mathbf{G}_{\varepsilon }\left( t\right) \mathbf{,}\boldsymbol{u}%
\left( t\right) \mathbf{,}\boldsymbol{v}\right) =\left( \boldsymbol{f}%
_{\varepsilon }\mathbf{\ }\left( t\right) \mathbf{,}\boldsymbol{v}\right) ,%
\end{array}%
\bigskip $

\noindent

\noindent This last relation being valid for all $m_{0},$
it remains true for all $\boldsymbol{v}\in \left\langle \boldsymbol{w}%
_{1,}\boldsymbol{w}_{2}\mathbf{,...,}\boldsymbol{w}_{m}\right\rangle ,$ $%
\forall m\in
%TCIMACRO{\U{2115} }%
%BeginExpansion
\mathbb{N}
%EndExpansion
^{\ast }.\medskip $

\smallskip

\noindent Finally let $\boldsymbol{v}\in V\mathbf{.}$ There exists $%
\boldsymbol{v}_{m}\in \left\langle \boldsymbol{w}_{1,}\boldsymbol{w}_{2}%
\mathbf{,...,}\boldsymbol{w}_{m}\right\rangle $ such that $\boldsymbol{v}%
_{m}\rightarrow \boldsymbol{v}$ in $V$ and

\begin{eqnarray}
&&\left( \boldsymbol{u}^{\prime }\left( t\right) ,\boldsymbol{v}\right) +\nu
\left( \left( \boldsymbol{u}\left( t\right) \mathbf{,}\boldsymbol{v}\right)
\right) +b\left( \boldsymbol{u}\left( t\right) \mathbf{,}\boldsymbol{u}%
\left( t\right) \mathbf{,}\boldsymbol{v}\right)  \notag \\
&&+b\left( \boldsymbol{u}\left( t\right) \mathbf{,\mathbf{G}}_{\varepsilon
}\left( t\right) \mathbf{,}\boldsymbol{v}\right) +b\left( \mathbf{G}%
_{\varepsilon }\left( t\right) \mathbf{,}\boldsymbol{u}\left( t\right)
\mathbf{,}\boldsymbol{v}\right) =\left( \boldsymbol{f}_{\varepsilon }\mathbf{%
\ }\left( t\right) \mathbf{,}\boldsymbol{v}\right)
\end{eqnarray}%
\bigskip

\noindent Now let us note that for all $t\in \left[ 0,T\right]
,\medskip $

$\ \ \ \ \ \ \ \ \ \ \ \ \ \ \ \ \ \ \ \ \ \ \ \ \ \ \ \ \ \ \ \ \ \ \ \ \ \
\ \ \ \
\begin{array}{lll}
\boldsymbol{u}_{m}\left( t\right) \rightarrow \boldsymbol{u}\left( t\right)
& \text{ weakly in } & V\text{,}%
\end{array}%
$

$\qquad \qquad $ \hspace{0.5cm}

\noindent and thus \medskip

$\ \ \ \ \ \ \ \ \ \ \ \ \ \ \ \ \ \ \ \ \ \ \ \ \ \ \ \ \ \ \
\begin{array}{lll}
\boldsymbol{u}_{m}\left( 0\right) =\boldsymbol{u}_{0m}\rightarrow
\boldsymbol{u}\left( 0\right) & \text{ weakly in } & V\text{.}%
\end{array}%
$

$\qquad \qquad $ \hspace{0.5cm}

\noindent Since

$\qquad \qquad \ \ \ \ \ \ \ \ \ \ \ \ \ \ \ \ \ \ \ \ \ \ \ \ \ \ \ \ \ \
\begin{array}{lll}
\boldsymbol{u}_{0m}\rightarrow \boldsymbol{u}_{0}^{\varepsilon }\hspace{0.5cm%
} & \text{in} & V\mathbf{,}%
\end{array}%
$

\noindent we have: $\ \ \ $

$\ \ \ \ \ \ \ \ \ \ \ \ \boldsymbol{u}\left( 0\right) =\boldsymbol{u}%
_{0}^{\varepsilon }.\bigskip $

\smallskip

\subsection{\textit{Existence of pressure.}}

\smallskip

\noindent From (31), one has, for all $\boldsymbol{v}\in
V\mathbf{,}$

\begin{equation*}
\left\langle \boldsymbol{u}^{\prime }-\nu \triangle \boldsymbol{u}+B\left(
\boldsymbol{u}\mathbf{,}\boldsymbol{u}\right) +B\left( \boldsymbol{u}\mathbf{%
,\mathbf{G}}_{\varepsilon }\right) +B\left( \mathbf{G}_{\varepsilon }\mathbf{%
,}\boldsymbol{u}\right) -\boldsymbol{f}_{\varepsilon }\mathbf{\ ,}\text{ }%
\boldsymbol{v}\right\rangle _{\mathbf{H}^{-1}(\Omega )\times \mathbf{H}%
_{0}^{1}(\Omega )}=0.
\end{equation*}

\smallskip

\noindent Consequently, there exists a unique function $p$ of
$L^{2}\left( 0,T\right) $ satisfying (16) and such that :$\qquad
\qquad $%
\begin{equation*}
p\in L^{2}\left( 0,T;\mathbf{H}^{1}\left( \Omega \right) \right) .
\end{equation*}

\bigskip

\noindent This ends the proof of theorem 2.1.$\square
\medskip $

\section{\protect\smallskip Uniqueness Theorem }

\noindent

\noindent \textbf{Theorem 3.1} \textit{Problem (16) has a unique
solution}.\medskip

\smallskip

\noindent \textbf{Proof.\medskip \smallskip }

\noindent Let $\boldsymbol{u}$ and $\boldsymbol{v}$ be two solutions
satisfying the hypotheses of theorem 2.1 and %
let $\boldsymbol{%
w}$ $\mathbf{=}$ $\boldsymbol{u}-\boldsymbol{v}\mathbf{.}$ Then one
has\medskip

$\qquad \dfrac{\partial \boldsymbol{w}}{\partial t}-\nu \triangle
\boldsymbol{w}$ $\mathbf{+}$ $\boldsymbol{w}\mathbf{.\nabla }\boldsymbol{u}$
$\mathbf{+}$ $\boldsymbol{v}\mathbf{.\nabla }\boldsymbol{w}$ $\mathbf{+}$ $%
\boldsymbol{w}\mathbf{.\nabla \mathbf{G}}_{\varepsilon }$ $\mathbf{+}$ $%
\mathbf{\mathbf{G}}_{\varepsilon }\mathbf{.\nabla }\boldsymbol{w}$ $\mathbf{%
=0\medskip }$

\noindent Multiplying by $\boldsymbol{w}\mathbf{,}$ we obtain
\medskip

$\qquad
\begin{array}{ll}
\dfrac{1}{2}\dfrac{d}{dt}\left\vert \boldsymbol{w}\left( t\right)
\right\vert ^{2}+\nu \left\Vert \boldsymbol{w}\left( t\right) \right\Vert
^{2}= & -\left( \boldsymbol{w}\mathbf{.\nabla }\boldsymbol{u}\mathbf{,}%
\boldsymbol{w}\right) -\left( \boldsymbol{v}\mathbf{.\nabla }\boldsymbol{w}%
\mathbf{\mathbf{,}}\boldsymbol{w}\right) \\
& -\left( \boldsymbol{w}\mathbf{.\nabla \mathbf{\mathbf{G}}_{\varepsilon },}%
\boldsymbol{w}\right) -\left( \mathbf{\mathbf{G}}_{\varepsilon }\mathbf{%
.\nabla }\boldsymbol{w}\mathbf{,}\boldsymbol{w}\right)%
\end{array}%
$

$\qquad $

\noindent But $b\left( \boldsymbol{v}\mathbf{,}\boldsymbol{w}\mathbf{%
\mathbf{,}}\boldsymbol{w}\right) =0$ and $b\left( \mathbf{\mathbf{G}}%
_{\varepsilon
}\mathbf{,}\boldsymbol{w}\mathbf{,}\boldsymbol{w}\right) =0.$ This
yields \medskip

$\qquad
\begin{array}{ll}
\dfrac{1}{2}\dfrac{d}{dt}\left\vert \boldsymbol{w}\left( t\right)
\right\vert ^{2}+\nu \left\Vert \boldsymbol{w}\left( t\right) \right\Vert
^{2}= & -b\left( \boldsymbol{w}\mathbf{,}\boldsymbol{u}\mathbf{,}\boldsymbol{%
w}\right) -b\left( \boldsymbol{w}\mathbf{,\mathbf{\mathbf{G}}_{\varepsilon },%
}\boldsymbol{w}\right) .%
\end{array}%
\hspace{1.5cm}$

$\qquad $

\noindent One then integrates with respect to t and we get \medskip

$\qquad \dfrac{1}{2}\left\vert \boldsymbol{w}\left( t\right)
\right\vert ^{2}+\nu \int_{0}^{t}\left\Vert \boldsymbol{w}\left(
s\right)
\right\Vert ^{2}ds=-\int_{0}^{t}b\left( \boldsymbol{w}\mathbf{,}%
\boldsymbol{u}\mathbf{,}\boldsymbol{w}\right) $ $ds-\int_{0}^{t}b%
\left( \boldsymbol{w}\mathbf{,\mathbf{\mathbf{G}}_{\varepsilon },}%
\boldsymbol{w}\right) $ $ds.\medskip $

\noindent Since \medskip

$\qquad
\begin{array}{lll}
\left\vert \int_{0}^{t}b\left( \boldsymbol{w}\mathbf{,}\boldsymbol{u}%
\mathbf{,}\boldsymbol{w}\right) ds\right\vert & \leq & C_{1}\int%
_{0}^{t}\left\Vert \boldsymbol{w}\left( s\right) \right\Vert _{%
\mathbf{L}^{4}\left( \Omega \right) }\left\Vert \boldsymbol{u}\left(
s\right) \right\Vert _{\mathbf{L}^{2}\left( \Omega \right) }ds \\
& \leq & C_{2}\int_{0}^{t}\left\vert \boldsymbol{w}\left( s\right)
\right\vert \text{ }\left\Vert \boldsymbol{w}\left( s\right)
\right\Vert
\left\Vert \boldsymbol{u}\left( s\right) \right\Vert ds \\
& \leq & \dfrac{\nu }{2}\int_{0}^{t}\left\Vert \boldsymbol{w}\left(
s\right) \right\Vert ^{2}ds+C_{3}\int_{0}^{t}\left\vert \boldsymbol{w%
}\left( s\right) \right\vert ^{2}\left\Vert \boldsymbol{u}\left( s\right)
\right\Vert ^{2}ds.%
\end{array}%
\medskip $

\noindent and, by the same way,  \medskip

$\qquad
\begin{array}{lll}
\int_{0}^{t}b\left( \boldsymbol{w}\mathbf{,\mathbf{\mathbf{G}}%
_{\varepsilon },}\boldsymbol{w}\right) ds & \leq & \dfrac{\nu }{2}%
\int_{0}^{t}\left\Vert \boldsymbol{w}\left( s\right) \right\Vert
^{2}ds+C_{4}\int_{0}^{t}\left\vert \boldsymbol{w}\left( s\right)
\right\vert ^{2}\left\vert \nabla \mathbf{\mathbf{\mathbf{G}}_{\varepsilon }}%
\left( s\right) \right\vert ^{2}ds.%
\end{array}%
\medskip $

\noindent it follows that \medskip

$\qquad \left\vert \boldsymbol{w}\left( t\right) \right\vert
^{2}\leq C_{5}\medskip \int_{0}^{t}\left\vert \boldsymbol{w}\left(
s\right)
\right\vert ^{2}\left( \left\vert \nabla \mathbf{\mathbf{\mathbf{G}}%
_{\varepsilon }}\left( s\right) \right\vert ^{2}+\left\Vert \boldsymbol{u}%
\left( s\right) \right\Vert ^{2}\right) ds\medskip $

\noindent Thanks to the Gronwall lemma, one deduces \medskip\ $%
\begin{array}{lll}
\boldsymbol{w} & = & 0.\square%
\end{array}%
\medskip  $

\section{Existence of strong reproductive solution \protect\medskip}

\noindent We first recall results obtained  by Kaniel et Shinbrot $%
\left[ 5\right] $ in the study of the following problem :%
\begin{equation}
\left\{
\begin{array}{lll}
\dfrac{\partial \boldsymbol{u}}{\partial t}-\nu \triangle \boldsymbol{u}%
\text{ }\mathbf{+}\text{ }\boldsymbol{u}\mathbf{.\nabla }\boldsymbol{u}\text{
}\mathbf{+}\text{ }\mathbf{\mathbf{\nabla }}p\text{ }\mathbf{=}\ \boldsymbol{%
f}\text{\ } & \text{in \ } & Q_{T} \\
\text{div }\boldsymbol{u}\text{ }\mathbf{=}\text{ }0 & \text{in} &
Q_{T}
\\
\boldsymbol{u}\text{ }\mathbf{=}\text{ }0 & \text{on} & \Sigma _{T} \\
\boldsymbol{u}\mathbf{(}0\mathbf{)=}\text{ }\boldsymbol{u}_{0} &
\text{in}
& \Omega%
\end{array}%
\right.
\end{equation}%
\medskip

\noindent where $\Omega $ is an open and bounded domain of $%
%TCIMACRO{\U{211d} }%
%BeginExpansion
\mathbb{R}
%EndExpansion
^{3},$ with a smooth boundary $\Gamma $. \medskip

\noindent The following result establishes the property of a
reproductive solution \medskip

\noindent \textbf{Theorem 4.1}. \textit{Let T $>$}$0,$%
\textit{  and }$\ \boldsymbol{f}$\ $\in \mathcal{B}_{R,T} $\ %
\textit{with }$\boldsymbol{f}$ \textit{small enough. }\ \textit{Then, %
there exists an unique function }$\boldsymbol{u}_{0}$,%
\textit{ independent %
of t, with }$\nabla \boldsymbol{u}_{0}\in \mathcal{B}_{R,T}$%
\textit{ and such that the solution of (32) reproduces its initial
value at }$ t=T$ :

\begin{equation*}
\boldsymbol{u}\left( \boldsymbol{x,}T\right) =\boldsymbol{u}%
\left( \boldsymbol{x,}0\right) =\boldsymbol{u}_{0}\left(
\boldsymbol{x}\right) ,
\end{equation*}%
\noindent\textit{where}
\begin{equation*}
\mathcal{B}_{R,T}=\left\{ \boldsymbol{u\in }\text{
}\mathbf{L}^{\infty }\left( 0,T;\mathbf{L}^{2}\left( \Omega \right)
\right) \text{ : }\left\Vert
\boldsymbol{u}\right\Vert _{\mathbf{L}^{\infty }\left( 0,T;\mathbf{L}%
^{2}\left( \Omega \right) \right) }\leq R\right\} .
\end{equation*}%
\medskip

\noindent We begin by recalling the following lemma.
\bigskip

\noindent \textbf{Lemma 4.2}.\textbf{\ }\textit{If }$\qquad $%
\begin{equation*}
\boldsymbol{u}\text{ }\mathbf{\in }\mathit{\ }\mathbf{L}^{2}\left( 0,T;%
\mathbf{H}^{2}\left( \Omega \right) \cap V\right) \mathit{\ and \ }\boldsymbol{%
u}^{\prime }\in \mathbf{L}^{2}\left( 0,T;\mathbf{H}\right)
\end{equation*}%
\textit{\ then }$\qquad $%
\begin{equation*}
\boldsymbol{u}\mathit{\ }\in C\left( \left[ 0,T\right] ;V\right)
\end{equation*}%
\textit{\ and }$\qquad $%
\begin{equation*}
\dfrac{d}{dt}\left\Vert \boldsymbol{u}\left( t\right) \right\Vert
^{2}=-2\left( \boldsymbol{u}^{\prime }\left( t\right) ,\triangle \boldsymbol{%
u}\left( t\right) \right) .\square
\end{equation*}

\bigskip

\noindent Now, let

\smallskip

\begin{equation}
\boldsymbol{v}_{0}\in \mathbf{H}^{1}(\Omega )\cap \mathbf{H,\quad }%
\boldsymbol{w}_{0}\text{ }\mathbf{\in \mathbf{H}}^{1}(\Omega )\text{ }%
\mathbf{\cap }\text{ }\mathbf{\mathbf{H,\quad }}\boldsymbol{g}\text{ }%
\mathbf{\in H}^{3/2,3/4}(\Sigma _{T})\text{ }
\end{equation}%
with

\begin{equation}
\text{ }\boldsymbol{g}.\boldsymbol{n}=0\text{ \thinspace on \quad
}\Sigma _{T}\text{ \quad and \quad }\boldsymbol{v}_{0}\left(
\boldsymbol{x}\right)
\text{ }\mathbf{=}\text{ }\boldsymbol{w}_{0}\left( \boldsymbol{x}\right) =%
\boldsymbol{g}\left( \boldsymbol{x},0\right) \text{ \quad
}\boldsymbol{x}\in \text{\thinspace }\Gamma .
\end{equation}

\bigskip

\noindent With these assumptions, it follows from theorem 2.1
that system (1), with data $\left( \boldsymbol{v}_{0},%
\boldsymbol{g}\right) ,$ (respectively $\left( \boldsymbol{w}_{0},%
\boldsymbol{g}\right) $), has an unique solution

\begin{equation*}
\boldsymbol{v}\text{ }\mathbf{\in }\mathit{\ }\mathbf{L}^{2}\left( 0,T;%
\mathbf{H}^{2}\left( \Omega \right) \cap \mathbf{H}\right) \cap \mathbf{L}%
^{\infty }\left( 0,T;\mathbf{H}^{1}\left( \Omega \right) \right) \text{ and }%
\boldsymbol{v}^{\prime }\in \mathbf{L}^{2}\left( 0,T;\mathbf{H}\right) ,
\end{equation*}

\bigskip

\noindent (respectively $\hspace{1.5cm}$%
\begin{equation*}
\boldsymbol{w}\text{ }\mathbf{\in }\mathit{\ }\mathbf{L}^{2}\left( 0,T;%
\mathbf{H}^{2}\left( \Omega \right) \cap \mathbf{H}\right) \cap \mathbf{L}%
^{\infty }\left( 0,T;\mathbf{H}^{1}\left( \Omega \right) \right)
\text{ and  \ }\boldsymbol{w}^{\prime }\in \mathbf{L}^{2}\left(
0,T;\mathbf{H}\right) \text{)}.
\end{equation*}

\bigskip

\noindent Let us now set $\mathbf{z=}$ $\boldsymbol{v}-\boldsymbol{w}%
\mathbf{.}$ Then

\begin{equation}
\left\{
\begin{array}{lll}
\dfrac{\partial \mathbf{z}}{\partial t}-\nu \triangle \mathbf{z+}\text{ }%
\boldsymbol{w}\mathbf{.\nabla z+z.\nabla }\boldsymbol{v}\text{ }\mathbf{+}%
\text{ }\mathbf{\mathbf{\nabla }}r\text{ }\mathbf{=0\quad } & \text{%
in \quad } & Q_{T}, \\
\text{div }\mathbf{z}=0 & \text{in} & Q_{T}, \\
\mathbf{z=}\text{ }0 & \text{on} & \text{ }\Sigma _{T}, \\
\mathbf{z(}0\mathbf{)=}\text{ }\boldsymbol{v}_{0}-\text{ }\boldsymbol{w}_{0}
& \text{in} & \text{ }\Omega .%
\end{array}%
\right.
\end{equation}

\bigskip

\noindent where $r=p-q$ \ ($q$ being the pressure corresponding to $%
\boldsymbol{w}$).\medskip \bigskip

\noindent \textbf{Lemma 4.3}. \textit{ If }

\begin{equation}
\max \left( \left\Vert \boldsymbol{v}\right\Vert _{\mathbf{L}^{\infty
}\left( 0,T;\mathbf{H}^{1}\left( \Omega \right) \right) },\left\Vert
\boldsymbol{w}\right\Vert _{\mathbf{L}^{\infty }\left( 0,T;\mathbf{H}%
^{1}\left( \Omega \right) \right) }\right) \leq M
\end{equation}%
\textit{under the assumptions (33) and (34) }\textit{with }$0<M<<1
$,\textit{\ then}

\begin{equation}
\dfrac{d}{dt}\left\Vert \mathbf{z}\left( t\right) \right\Vert
^{2}+\nu \left\Vert \mathbf{z}\left( t\right) \right\Vert ^{2}\leq 0
\end{equation}%
\textit{and thus, for all }$t\in \left[ 0,T\right],$
\begin{equation}
\left\Vert \boldsymbol{v}\left( t\right) -\boldsymbol{w}\left(
t\right) \right\Vert \leq \left\Vert
\boldsymbol{v}_{0}-\boldsymbol{w}_{0}\right\Vert \exp \left( -\nu
t\right) .
\end{equation}

\smallskip

\noindent \textbf{Proof.}\bigskip

\noindent Let P: $\mathbf{L}^{2}\left( \Omega \right) \rightarrow \mathbf{H,%
}$ be the orthogonal projection operator.\ Then \medskip

$\qquad \qquad \forall \mathbf{\varphi \in H,}\left( \mathbf{\mathbf{\nabla }%
}r\text{, }\mathbf{\varphi }\right) =0.$

\smallskip

\noindent In particular, let us multiply (35) by $P\triangle \mathbf{z=}$ $A%
\mathbf{z:\medskip }$

$\qquad \qquad \dfrac{1}{2}\dfrac{d}{dt}\left\Vert \mathbf{z}\left( t\right)
\right\Vert ^{2}+\nu \left\vert A\mathbf{z}\right\vert ^{2}=-\left(
\boldsymbol{w}\mathbf{.\nabla z,}A\mathbf{z}\right) -\left( \mathbf{z.\nabla
}\boldsymbol{v}\mathbf{\mathbf{,}}A\mathbf{z}\right) $

\noindent But \medskip

$\qquad
\begin{array}{ll}
\left\vert \left( \boldsymbol{w}\mathbf{.\nabla z,}A\mathbf{z}\right)
\right\vert & \leq \left\Vert \boldsymbol{w}\right\Vert _{\mathbf{L}%
^{4}\left( \Omega \right) }\left\Vert \nabla \mathbf{z}\right\Vert _{\mathbf{%
L}^{4}\left( \Omega \right) }\left\vert A\mathbf{z}\right\vert \\
& \leq C\left\Vert \boldsymbol{w}\right\Vert \left\vert A\mathbf{z}%
\right\vert ^{2}%
\end{array}%
$

$\qquad \qquad \qquad \qquad $\ $\ \ \ \ \ \ \ \ \ \ \ \ \ \ \ \ \ $

\noindent and \medskip

$\qquad
\begin{array}{ll}
\left\vert \left( \mathbf{z.\nabla }\boldsymbol{v}\mathbf{\mathbf{,}}A%
\mathbf{z}\right) \right\vert & \leq \left\Vert \mathbf{z}\right\Vert _{%
\mathbf{L}^{\infty }\left( \Omega \right) }\left\Vert \boldsymbol{v}%
\right\Vert _{{}}\left\vert A\mathbf{z}\right\vert \\
& \leq C\left\Vert \boldsymbol{v}\right\Vert _{{}}\left\vert A\mathbf{z}%
\right\vert ^{2}.%
\end{array}%
$

$\qquad \qquad $

\noindent So that if \medskip

$\qquad \qquad C\left( \left\Vert \boldsymbol{v}\right\Vert _{\mathbf{L}%
^{\infty }\left( 0,T;\mathbf{H}^{1}\left( \Omega \right) \right)
}+\left\Vert \boldsymbol{w}\right\Vert _{\mathbf{L}^{\infty }\left( 0,T;%
\mathbf{H}^{1}\left( \Omega \right) \right) }\right) \leq \dfrac{\nu }{2}$

\noindent then \medskip

$\qquad \qquad \dfrac{d}{dt}\left\Vert \mathbf{z}\left( t\right) \right\Vert
^{2}+\nu \left\Vert \mathbf{z}\left( t\right) \right\Vert ^{2}\leq 0\medskip
$

\noindent and one deduces (38). $\ \square \medskip
\medskip
\medskip $
\subsection{The main result}

\noindent \textbf{Lemma 4.4}. \textit{Suppose that }\textbf{g }\textit{and }\textbf{v}%
$_{0}$ \textit{satisfy hypotheses (4)-(5) and (9).\ Let us suppose
moreover that }\textbf{f}$_{\varepsilon }\in \mathbf{L}^{\infty
}\left( 0,T;\mathbf{L}^{2}\left( \Omega \right) \right) $
\textit{and that }
\begin{equation}
\left\Vert \boldsymbol{g}\right\Vert _{\mathbf{H}^{3/2,3/4}(\Sigma
_{T})}+\left\Vert \boldsymbol{v}_{0}\right\Vert
_{\mathbf{H}^{1}(\Omega )}\leq \alpha
\end{equation}%
%\medskip
\begin{equation}
\left\Vert \boldsymbol{f}_{\varepsilon }\right\Vert
_{\mathbf{L}^{\infty }\left( 0,T;\mathbf{L}^{2}\left( \Omega \right)
\right) }\leq K
\end{equation}%
%\medskip
\noindent \textit{with} $\alpha >0$ and $0<K<<1 $ \textit{}.\
\textit{Then, if }$\boldsymbol{u}$\textit{\ is the solution given by
theorem 2.1, one has:}

\begin{equation}
\underset{t\in \left[ 0,T\right] }{\sup }\left\Vert \nabla \boldsymbol{u}%
\left( t\right) \right\Vert _{\mathbf{L}^{2}\left( \Omega \right) }\leq M
\end{equation}%
\bigskip

\noindent \textbf{Remark 4.5. } Let us recall that \medskip

\begin{equation*}
\boldsymbol{u}_{0}=\boldsymbol{v}_{0}-\mathbf{G}_{\varepsilon }\left(
.,0\right)
\end{equation*}%
\bigskip Consequently, if hypothesis (39) takes place, one has from lemma 1.4 :\medskip

$\qquad
\begin{array}{ll}
\ \left\Vert \boldsymbol{u}_{0}\right\Vert & \leq \left\Vert \boldsymbol{u}%
_{0}\right\Vert _{\mathbf{H}^{1}(\Omega )}\leq \left\Vert \boldsymbol{v}%
_{0}\right\Vert _{\mathbf{H}^{1}(\Omega )}+\left\Vert \mathbf{G}%
_{\varepsilon }\left( .,0\right) \right\Vert _{\mathbf{H}^{1}(\Omega )} \\
& \leq \left\Vert \boldsymbol{v}_{0}\right\Vert _{\mathbf{H}^{1}(\Omega
)}+L\left( \left\Vert \boldsymbol{g}\right\Vert _{\mathbf{H}%
^{3/2,3/4}(\Sigma _{T})}+\left\Vert \boldsymbol{v}_{0}\right\Vert _{\mathbf{H%
}^{1}(\Omega )}\right) \\
& \leq \alpha \left( L+1\right) =M.\square%
\end{array}%
\bigskip $

\noindent \textbf{Proof of lemma 4.4.\bigskip } \textrm{(see Batchi
$ \left[ 5\right] $)}

\noindent Let us multiply (16) by A\textbf{\textit{u}} and integrate
on $\Omega $ :\medskip

$\qquad
\begin{array}{ll}
\ \dfrac{1}{2}\dfrac{d}{dt}\left\Vert \boldsymbol{u}\right\Vert ^{2}+\nu
\left\vert A\boldsymbol{u}\right\vert ^{2}\leq & \ \ \int_{\Omega }%
\boldsymbol{f}_{\varepsilon }\mathbf{\ }.A\boldsymbol{u}dx-\int_{%
\Omega }\left( \boldsymbol{u}\mathbf{.}\nabla \boldsymbol{u}\right) .A%
\boldsymbol{u}dx \\
& -\int_{\Omega }\left( \boldsymbol{u}\mathbf{.}\nabla \mathbf{G}%
_{\varepsilon }\right) .A\boldsymbol{u}dx-\int_{\Omega }\left(
\mathbf{G}_{\varepsilon }\mathbf{.}\nabla \boldsymbol{u}\right) .A%
\boldsymbol{u}dx%
\end{array}%
$

$\ \ $

\noindent But \medskip

$\qquad
\begin{array}{ll}
\ \ \ \left\vert \int_{\Omega }\left( \boldsymbol{u}\mathbf{.}\nabla
\boldsymbol{u}\right) .A\boldsymbol{u}\text{ }dx\right\vert & \leq
\left\Vert \boldsymbol{u}\right\Vert _{\mathbf{L}^{\infty }\left(
\Omega
\right) }\left\Vert \boldsymbol{u}\right\Vert \left\vert A\boldsymbol{u}%
\right\vert \\
& \leq C_{1}\left\Vert \boldsymbol{u}\right\Vert \left\vert A\boldsymbol{u}%
\right\vert ^{2},%
\end{array}%
$

$\ \ \ \ \ \ \ \ \ \ \ \ \ \ \ \ \ \ \ \ \ \ \ \ \ \ \ \ \ \ \ \ \ $

\noindent where $C_{1}$ is such that $\left\Vert
\boldsymbol{u}\right\Vert
_{\mathbf{L}^{\infty }\left( \Omega \right) }\leq $ $C_{1}\left\vert A%
\boldsymbol{u}\right\vert .\medskip $

\noindent In the same way, one also has \medskip

$\ \ \ \qquad \left\vert \int_{\Omega }\left( \boldsymbol{u}\mathbf{.%
}\nabla \mathbf{G}_{\varepsilon }\right) .A\boldsymbol{u}dx\right\vert \leq
C_{1}\left\Vert \nabla \mathbf{G}_{\varepsilon }\right\Vert _{\mathbf{L}%
^{2}\left( \Omega \right) }$ $\left\vert A\boldsymbol{u}\right\vert
^{2}\medskip $

\noindent But thanks to the lemma 1.4, one knows that \medskip

\ \ \ \textbf{G}$_{\varepsilon }\in \mathbf{L}^{\infty }\left( 0,T;\mathbf{H}%
^{1}\left( \Omega \right) \right) $

\noindent and \medskip

$\qquad
\begin{array}{ll}
\ \left\Vert \nabla \mathbf{G}_{\varepsilon }\right\Vert _{\mathbf{L}%
^{2}\left( \Omega \right) } & \leq C_{2}\left\Vert \mathbf{G}_{\varepsilon
}\right\Vert _{\mathbf{H}^{2,1}(Q_{T})} \\
& \leq C_{2}L\left( \left\Vert \boldsymbol{g}\right\Vert _{\mathbf{H}%
^{3/2,3/4}(\Sigma _{T})}+\left\Vert \boldsymbol{v}_{0}\right\Vert _{\mathbf{H%
}^{1}(\Omega )}\right) \\
& \leq C_{3}\alpha .%
\end{array}%
$

$\ \ \ \ \ \ \ \ \ \ \ \ \ \ \ \ \ \ \ \ \ \ \ \ \ \ $

\noindent It then follows that \medskip

$\qquad
\begin{array}{ll}
\ \left\vert \int_{\Omega }\left( \mathbf{G}_{\varepsilon }\mathbf{.}%
\nabla \boldsymbol{u}\right) .A\boldsymbol{u}dx\right\vert & \leq \left\Vert
\mathbf{G}_{\varepsilon }\right\Vert _{\mathbf{L}^{4}\left( \Omega \right)
}\left\Vert \nabla \boldsymbol{u}\right\Vert _{\mathbf{L}^{4}\left( \Omega
\right) }\left\vert A\boldsymbol{u}\right\vert \\
& \leq C_{4}\left\Vert \mathbf{G}_{\varepsilon }\right\Vert _{\mathbf{H}%
^{1}\left( \Omega \right) }\left\vert A\boldsymbol{u}\right\vert \left\Vert
\nabla \boldsymbol{u}\right\Vert _{\mathbf{L}^{2}\left( \Omega \right)
}^{1/2}\left\Vert \nabla ^{2}\boldsymbol{u}\right\Vert _{\mathbf{L}%
^{2}\left( \Omega \right) }^{1/2} \\
& \leq C_{5}\alpha \left\Vert \boldsymbol{u}\right\Vert
_{{}}^{1/2}\left\vert A\boldsymbol{u}\right\vert ^{3/2} \\
& \leq C_{5}\alpha \sqrt{C_{6}}\left\vert A\boldsymbol{u}\right\vert ^{2},%
\end{array}%
$

$\ \ \ \ \ \ \ \ \ \ \ \ \ \ \ \ \ \ \ \ \ \ \ \ \ \ \ \ \ \ \ \ \ $

\noindent with $\ \ \left\Vert \boldsymbol{u}\right\Vert \leq
C_{6}\left\vert A\boldsymbol{u}\right\vert .\medskip $

\noindent Thus,
\begin{equation}
\dfrac{1}{2}\dfrac{d}{dt}\left\Vert \boldsymbol{u}\right\Vert ^{2}+\nu
\left\vert A\boldsymbol{u}\right\vert ^{2}\leq \left\vert \boldsymbol{f}%
_{\varepsilon }\right\vert \text{ }\left\vert A\boldsymbol{u}\right\vert
+C_{1}\left\Vert \boldsymbol{u}\right\Vert \text{ }\left\vert A\boldsymbol{u}%
\right\vert ^{2}+C_{1}C_{3}\alpha \left\vert A\boldsymbol{u}\right\vert
^{2}+C_{5}\alpha \sqrt{C_{6}}\left\vert A\boldsymbol{u}\right\vert ^{2}.
\end{equation}%
\medskip \bigskip Let $\ \varphi \left( t\right) =\left\Vert \boldsymbol{u%
}\left( t\right) \right\Vert \medskip $

\noindent i) Let us first suppose that $\left\Vert
\boldsymbol{u}_{0}\right\Vert <M.\bigskip $

\noindent Let $t_{0}>0$ be the smallest $t>0$ \ such that $\varphi
\left( t_{0}\right) =M.$ According to (41), one then has \medskip

$\qquad
\begin{array}{ll}
\dfrac{1}{2}\dfrac{d}{dt}\left\Vert \boldsymbol{u}\left( t\right)
\right\Vert _{t=t_{0}}^{2}+\nu \left\vert A\boldsymbol{u}\left( t_{0}\right)
\right\vert ^{2}\leq & K\left\vert A\boldsymbol{u}\left( t_{0}\right)
\right\vert +C_{1}M\left\vert A\boldsymbol{u}\left( t_{0}\right) \right\vert
^{2} \\
& +C_{1}C_{3}\alpha \left\vert A\boldsymbol{u}\left( t_{0}\right)
\right\vert ^{2}+C_{5}\alpha \sqrt{C_{6}}\left\vert A\boldsymbol{u}\left(
t_{0}\right) \right\vert ^{2}.%
\end{array}%
$

$\ \ \ \ \ \ \ \ $

\noindent Let us choose $\alpha $ \ sufficiently small and $K$ such
that
\medskip

$\ \ \ \ \ \ K=\dfrac{\nu }{8}\dfrac{1}{C_{6}}M,\qquad \left(
C_{1}M+C_{1}C_{3}\alpha +C_{5}\alpha \sqrt{C_{6}}\right) \leq \dfrac{3\nu }{8%
}$

\noindent Then \medskip

$\qquad
\begin{array}{ll}
\ \ \dfrac{1}{2}\dfrac{d}{dt}\left\Vert \boldsymbol{u}\left( t\right)
\right\Vert _{t=t_{0}}^{2}+\nu \left\vert A\boldsymbol{u}\left( t_{0}\right)
\right\vert ^{2} & \leq \dfrac{\nu }{8}\dfrac{1}{C_{6}}M\left\vert A%
\boldsymbol{u}\left( t_{0}\right) \right\vert +\dfrac{3\nu }{8}\left\vert A%
\boldsymbol{u}\left( t_{0}\right) \right\vert ^{2} \\
\dfrac{1}{2}\dfrac{d}{dt}\left\Vert \boldsymbol{u}\left( t\right)
\right\Vert _{t=t_{0}}^{2}+\nu \left\vert A\boldsymbol{u}\left( t_{0}\right)
\right\vert ^{2} & \leq \dfrac{\nu }{8}\dfrac{1}{C_{6}}\left\Vert
\boldsymbol{u}\left( t_{0}\right) \right\Vert \left\vert A\boldsymbol{u}%
\left( t_{0}\right) \right\vert +\dfrac{3\nu }{8}\left\vert A\boldsymbol{u}%
\left( t_{0}\right) \right\vert ^{2} \\
\dfrac{1}{2}\dfrac{d}{dt}\left\Vert \boldsymbol{u}\left( t\right)
\right\Vert _{t=t_{0}}^{2}+\nu \left\vert A\boldsymbol{u}\left( t_{0}\right)
\right\vert ^{2} & \leq \dfrac{\nu }{2}\left\vert A\boldsymbol{u}\left(
t_{0}\right) \right\vert ^{2}.%
\end{array}%
$

$\ \ \ \ \ \ \ \ \ \ \ \ \ \ \ \ \ \ \ \ \ \ \ \ \ \ \ \ \ \ \ \ \ \ \ \ \ \
\ \ \ \ \ \ \ \ \ \ \ \ \ \ \ $

\noindent Thus \medskip

$\ \ \ \ \ \ \ \ \dfrac{d}{dt}\left\Vert \boldsymbol{u}\left( t\right)
\right\Vert _{t=t_{0}}^{2}+\nu \left\vert A\boldsymbol{u}\left( t_{0}\right)
\right\vert ^{2}\leq 0\medskip $

\noindent which implies that \medskip

$\ \ \ \ \ \ \ \dfrac{d}{dt}\left\Vert \boldsymbol{u}\left( t\right)
\right\Vert _{t=t_{0}}^{2}\leq 0\medskip $

\noindent Consequently, there exists $t^{\ast }\in \left[
0,t_{0}\right[ $ such that \medskip

$\ \ \ \ \ \ \ \ \varphi \left( t^{\ast }\right) >\varphi \left(
t_{0}\right) ,$ in contradiction with the definition of $%
t_{0}.\medskip \ $

\noindent Therefore \medskip

$\ \ \ \ \ \ \ \ \forall t\in \left[ 0,T\right] ,$ $\ \varphi \left(
t\right) <M.\medskip \bigskip $

\noindent ii) Suppose now that $\left\Vert \boldsymbol{u}%
_{0}\right\Vert =M.\bigskip $

\noindent According to the above calculations, one verifies that $%
\varphi ^{\prime }\left( 0\right) <0$ and thus there exists $t^{\ast
}>0$ such that

$\ \ \ \ \ \ \ \ $%
\begin{equation*}
\forall t\in \left] 0,t^{\ast }\right] ,\varphi \left( t\right) <M.
\end{equation*}

\bigskip

\noindent Repeating the reasoning made in i)$,$ one shows that on $%
\left[ t^{\ast },T\right] $, $\varphi \left( t\right) <M,$ and this
ends the proof.$\square
\bigskip
\smallskip $

\noindent \textbf{Remark 4.6. }From now on, we assume that
$\boldsymbol{g}$\textbf{\ }\ does not dependent on time.\ More
precisely, it is supposed that
\begin{equation}
\boldsymbol{g}\mathbf{\in  H }^{3/2}\left( \Gamma \right) ,\text{ \ \ \ }%
\boldsymbol{g}\mathbf{.}\boldsymbol{n}\mathbf{=}\text{ }0\text{ on
}\Gamma .
\end{equation}

\medskip

\noindent One recalls that $\boldsymbol{v}_{0}\in
\mathbf{H}^{1}\left( \Omega \right) $ satisfies

\begin{equation}
\text{div }\boldsymbol{v}_{0}=0\text{ in }\Omega ,\text{ \ \ }\boldsymbol{v%
}_{0}.\boldsymbol{n}\mathbf{=}\text{ }0\text{ on }\Gamma
\end{equation}%
and that

\begin{equation}
\text{\ }\boldsymbol{v}_{0}\text{ }\mathbf{=\ }\boldsymbol{g}\text{
\ \ \ on }\Gamma .
\end{equation}

\bigskip

\noindent One knows that there exists \ $\mathbf{G\in H}^{2}\left(
\Omega \right) $ \ such that

\begin{equation}
\left\{
\begin{array}{c}
\text{div }\mathbf{G}=0\text{ \ \ \ \ \ \ in }\Omega , \\
\text{\ }\mathbf{G}\text{ }\mathbf{=\ }\boldsymbol{g}\text{ \ \ \ \ \ \ \ \
on \ }\Gamma ,%
\end{array}%
\right.
\end{equation}%
with

\begin{equation}
\left\Vert \mathbf{G}\right\Vert _{\mathbf{H}^{2}\left( \Omega
\right) }\leq C\left\Vert \boldsymbol{g}\right\Vert
_{\mathbf{H}^{3/2}\left( \Gamma \right) }.
\end{equation}

\bigskip

\noindent Processing as in lemma 1.4, one shows the existence, for
all $\varepsilon >0,$ of $\mathbf{G}_{\varepsilon }\in
\mathbf{H}^{2}\left( \Omega \right) $ satisfying (44)-(47) and the
estimates:

\begin{equation}
\forall \boldsymbol{v}\text{ }\mathbf{\in \boldsymbol{V},}\text{ }\left\vert
b\left( \boldsymbol{v},\mathbf{G}_{\varepsilon },\boldsymbol{v}\right)
\right\vert \leq \varepsilon \left\Vert \boldsymbol{g}\right\Vert ^{2}
\end{equation}

\bigskip

\noindent The right side $\boldsymbol{f}_{\varepsilon }$ in system
(16) then becomes independent of time and satisfies

\begin{equation}
\boldsymbol{f}_{\varepsilon }\in L^{\infty }\left( 0,T;L^{2}\left(
\Omega \right) ^{2}\right)
\end{equation}%
In the same way, \ $\boldsymbol{u}_{0}^{\varepsilon }$ becomes

\begin{equation}
\ \boldsymbol{u}_{0}^{\varepsilon }=\text{\ }\boldsymbol{v}_{0}-\mathbf{G}%
_{\varepsilon }
\end{equation}

\medskip

\noindent with $\mathbf{G}_{\varepsilon }$ depends only on $%
\boldsymbol{g}\mathbf{.\square }$

\medskip
\subsection{Reproductive solution result}
\noindent With these assumptions on $\boldsymbol{g}$ and $\boldsymbol{v}%
_{0} $, lemma 4.2 remains naturally valid and one is able to
establish the theorem which follows :\bigskip

\noindent \textbf{Theorem 4.7. }\textit{Let }$\boldsymbol{g}%
\mathbf{\in H}^{3/2}(\Gamma )$ \textit{such that }$\boldsymbol{g}\mathbf{.}%
\boldsymbol{n}\mathbf{=}$ $0$ \textit{on } $\Gamma $ \textit{ and }

\begin{equation}
\left\Vert \boldsymbol{g}\right\Vert _{\mathbf{H}^{3/2}\left( \Gamma
\right) }\leq \alpha \text{ }
\end{equation}%
%\smallskip
\noindent \textit{with }$ 0< \alpha <<1 $ \textit{.
Then, there exists } $\boldsymbol{v}_{0}\in \mathbf{H}^{1}\left( \Omega \right) $ \textit{%
such that} div $\boldsymbol{v}_{0}=0$ \ \ \textit{in} $\Omega $
\textit{and
}$\ \boldsymbol{v}_{0}$ $\mathbf{=\ }\boldsymbol{g}$ \ \ \textit{on} $%
\Gamma ,$ \textit{and such that the solution }$\ $\textit{\ }$\boldsymbol{v}$ $%
\mathbf{=\ }\boldsymbol{u}$ $\mathbf{+}$ $\mathbf{G}_{\varepsilon }$\textit{%
\ where }$\boldsymbol{u}$ \textit{is} \textit{given by theorem %
2.1, is reproductive:}
\begin{equation*}
\ \boldsymbol{v}\left( T\right) =\ \boldsymbol{v}\left( 0\right) =\ \
\boldsymbol{v}_{0}.
\end{equation*}%
\medskip
%\smallskip
\noindent \textbf{Proof. } \ Let $\mathbf{G}_{\varepsilon }\in
\mathbf{H}^{2}\left( \Omega \right) $ be the extension of $\boldsymbol{g}$ satisfying%
(45)-(47) and

\begin{equation*}
\boldsymbol{f}_{\varepsilon }\text{ }\mathbf{=}\text{ }\nu \triangle \mathbf{%
G}_{\varepsilon }-\mathbf{\mathbf{G}}_{\varepsilon }\mathbf{.\nabla \mathbf{G%
}}_{\varepsilon }
\end{equation*}

\bigskip

\noindent Let $\boldsymbol{u}_{0}^{\varepsilon }=$\ $\boldsymbol{v}_{0}-%
\mathbf{G}_{\varepsilon }\in V$ and $\ \boldsymbol{u}$ $\in $
$L^{2}\left( 0,T;\mathbf{H}^{2}\left( \Omega \right) \right) \cap
L^{\infty }\left(
0,T;V\right) $ be the unique solution of (16).\ We note that the function $%
\boldsymbol{v}$ $\mathbf{=\ }\boldsymbol{u}$ $\mathbf{+}$ $\mathbf{G}%
_{\varepsilon }$ is the unique solution of the initial problem (1).
As in the proof of lemma 4.3, it is clear that if $\left\Vert
\boldsymbol{u}_{0}^{\varepsilon }\right\Vert <M,$ then

\begin{equation*}
\underset{t\in \left[ 0,T\right] }{\sup }\left\Vert \boldsymbol{u}\left(
t\right) \right\Vert \leq M
\end{equation*}

\bigskip

\noindent provided that $\left\Vert \boldsymbol{f}_{\varepsilon }\right\Vert _{%
\mathbf{L}^{2}\left( \Omega \right) }$ is sufficiently small, which
follows from (49).\bigskip

\noindent Let us define the application

\begin{equation*}
\begin{array}{ll}
\text{L : } & \text{\ }\boldsymbol{u}_{0}^{\varepsilon }\longrightarrow
\boldsymbol{u}\left( .,T\right) \\
& B_{M}\text{ }\longrightarrow \text{ }B_{M}%
\end{array}%
\end{equation*}%
where $\ \ \ \ B_{M}=\left\{ \mathbf{z\in }\boldsymbol{V}\mathbf{,}\text{ }%
\left\Vert \mathbf{z}\right\Vert \leq M\right\} ;\medskip $

\noindent $\boldsymbol{u}\left( .,T\right) $ being the unique
solution of (16) at $t=T.$\smallskip

\noindent Moreover, as in remark 4.5, it is clear that if $%
\left\Vert \boldsymbol{v}_{0}\right\Vert \leq \alpha $ and
$\left\Vert
\boldsymbol{w}_{0}\right\Vert \leq \alpha $ then%
\begin{equation*}
\begin{array}{lll}
\left\Vert \boldsymbol{u}_{0}^{\varepsilon }\right\Vert \leq M &
\text{and} &
\left\Vert \boldsymbol{w}_{0}^{\varepsilon }\right\Vert \leq M,%
\end{array}%
\end{equation*}

\medskip

\noindent with $\ \ \ \mathbf{y}_{0}^{\varepsilon }=$\ $\boldsymbol{w}_{0}-%
\mathbf{G}_{\varepsilon }.\medskip $

\noindent So that%
\begin{equation*}
\begin{array}{ll}
\text{L}\boldsymbol{u}_{0}^{\varepsilon }\left( t\right) -\text{L}\mathbf{y}%
_{0}^{\varepsilon }\left( t\right) & =\boldsymbol{u}\left( t\right) -\mathbf{%
y}\left( t\right) \\
& =\boldsymbol{u}\left( t\right) -\mathbf{G}_{\varepsilon }-\left( \mathbf{y}%
\left( t\right) -\mathbf{G}_{\varepsilon }\right) \\
& =\boldsymbol{v}\left( t\right) -\boldsymbol{w}\left( t\right) ,%
\end{array}%
\end{equation*}

\bigskip

\noindent and, according to lemma 4.2%
\begin{equation*}
\begin{array}{ll}
\left\Vert \text{L}\boldsymbol{u}_{0}^{\varepsilon }\left( t\right) -\text{L}%
\mathbf{y}_{0}^{\varepsilon }\left( t\right) \right\Vert & =\left\Vert
\boldsymbol{v}\left( T\right) -\boldsymbol{w}\left( T\right) \right\Vert \\
& \leq \left\Vert \ \boldsymbol{v}_{0}-\ \boldsymbol{w}_{0}\right\Vert \exp
\left( -\nu T\right) \\
& \leq \left\Vert \boldsymbol{u}_{0}^{\varepsilon }-\mathbf{y}%
_{0}^{\varepsilon }\right\Vert \exp \left( -\nu T\right)%
\end{array}%
\end{equation*}%
Thus L is a contraction and has a fixed point.$\square
\medskip $

\end{document}